\documentclass[11pt]{article}

\usepackage{amsfonts,latexsym,amstext}
\usepackage{amsmath,amscd,amssymb}

\newtheorem{theorem}{Theorem}[section]
\newtheorem{lemma}[theorem]{Lemma}
\newtheorem{proposition}[theorem]{Proposition}
\newtheorem{definition}{Definition}[section]

\newtheorem{hypothesis}[theorem]{Hypothesis}

\newtheorem{remark}[theorem]{Remark}
\newtheorem{corollary}[theorem]{Corollary}
\newenvironment{proof}[1][Proof]{\medskip \noindent \textbf{#1.} }{\ \rule{0.5em}{0.5em}}
\makeatletter \@addtoreset{equation}{section}

\makeatother

\def\sqr#1#2{{\vcenter{\vbox{\hrule height .#2pt \hbox{\vrule
 width .#2pt height#1pt \kern#1pt \vrule
width .#2pt} \hrule height .#2pt}}}}

\def\ds{\begin{displaystyle}}
\def\eds{\end{displaystyle}}
\def\dis{\displaystyle }
\def\<{\langle }
\def\>{\rangle }

\def\R{\mathbb R}

\def\E{\mathbb E}
\def\P{\mathbb P}

\def\calf{{\cal F}}
\def\calg{{\cal G}}

\title{Ergodic BSDEs and Optimal Ergodic Control in Banach Spaces}

\date{}
\author{Marco Fuhrman
\\
Dipartimento di Matematica,
Politecnico di Milano\\
piazza Leonardo da Vinci 32, 20133 Milano, Italy\\
e-mail: marco.fuhrman@polimi.it
\\ \\
Ying Hu\\
IRMAR, Universit\'e Rennes 1\\ Campus de Beaulieu, 35042 RENNES
Cedex, France\\ e-mail: ying.hu@univ-rennes1.fr
\\ \\
Gianmario Tessitore
\\
Dipartimento di Matematica e Applicazioni, Universit\`a di Milano-Bicocca\\
Via Cozzi 53, 20135, Milano Italy\\
e-mail: gianmario.tessitore@unimib.it }

\begin{document}

\maketitle

\begin{abstract}

In this paper we introduce a new kind of Backward
Stochastic Differential Equations, called ergodic
BSDEs, which arise naturally in the study of optimal ergodic
control. We study the existence, uniqueness and regularity of
solution to
 ergodic
BSDEs. Then we apply these results to the optimal ergodic control
of a Banach valued stochastic state equation. We also establish
the link between the ergodic BSDEs and the associated
Hamilton-Jacobi-Bellman equation. Applications are given
to ergodic control of stochastic partial differential equations.

\end{abstract}

\section{Introduction}

In this paper we study the following type of (markovian) backward
stochastic differential equations with infinite horizon (that we
shall call \textit{ergodic} BSDEs or EBSDEs for short):
\begin{equation}\label{EBSDE*}
Y^x_t=Y^x_T +\int_t^T\left[\psi(X^x_\sigma,Z^x_\sigma)-
\lambda\right]d\sigma-\int_T^T Z^x_\sigma dW_\sigma, \quad 0\le
t\le T <\infty.
\end{equation}
In equation (\ref{EBSDE*}) $X^x$ is the solution of a forward
stochastic differential equation with
 values in a Banach space $E$
starting at $x$ and $(W_t)_{t\geq 0}$ is a cylindrical Wiener
process in a Hilbert space $\Xi$.

Our aim is to find a triple
 $(Y,Z,\lambda)$, where $Y,Z$ are
 adapted processes taking values in
 $\mathbb{R}$ and $\Xi^*$ respectively and $\lambda$
  is a real number. $\psi:E\times \Xi^*\to \R$ is
  a given function.
We stress the fact that $\lambda$ is part of the unknowns of
equation (\ref{EBSDE*}) and this is the reason  why the above is a
new class of BSDEs.

$ $

It is by now well known that BSDEs provide an efficient
alternative tool to study optimal control problems, see, e.g.
\cite{peng93}, \cite{ElKaMaz} or, in an infinite dimensional
framework, \cite{{FuTe1}}, \cite{masiero}. But up to our best
knowledge, there exists no work in which  BSDE techniques are
applied to
 optimal control problems with \emph{ergodic} cost functionals that is
 functionals depending only on the asymptotic behavior of the state
 (see e.g. the cost defined in formula
 (\ref{ergodic-cost*}) below).

$ $

\noindent The purpose of the present paper is to show that
backward
 stochastic differential equations, in particular
 the class of EBSDEs mentioned above,  are a very useful tool in the
 treatment of ergodic control problems as
 well, especially in an infinite dimensional framework.

 $ $

 \noindent There is a fairly large amount of literature dealing by
 analytic techniques with optimal ergodic control problems
  for finite dimensional stochastic state
 equations.
 We just mention the basic papers by Bensoussan and Frehse
  \cite{BeFr} and by Arisawa and Lions \cite{ArLi} where the
problem is treated through the study of  the corresponding
Hamilton-Jacobi-Bellman (HJB) equation (solutions  are understood
 in a classical sense and in a viscosity sense,
respectively).

Concerning the infinite dimensional case it is known that both
classical and viscosity notions of solutions are not so suitable
concepts. Maslowski and Goldys in \cite{GoMa} employ a mild
formulation of the Hamilton-Jacobi-Bellman equation in a
 Hilbertian framework (see \cite{C} and references within for the
corresponding mild formulations  in the standard cases). In
\cite{GoMa} the authors prove, by a fixed point argument that
exploits the smoothing properties of the Ornstein-Uhlenbeck
semigroup corresponding to the state equation,  existence and
uniqueness of the  solution of the stationary HJB equation for
discounted infinite horizon costs. Then they pass to the limit, as
the discount goes to zero, to obtain a mild solution of the  HJB
equation for the ergodic problem (see also \cite{duncans}). Such
techniques need to assume, beside natural condition on  the
dissipativity of the state equation, also non-degeneracy of the
noise and a limitation on the lipschitz constant (with respect to
the gradient variable)  of the hamiltonian function. This last
condition carries a bound on the size of the control domain (see
\cite{FuTe-ell} for similar conditions in the infinite horizon
case).

$ $

The introduction of EBSDEs allow us to treat Banach valued state
equations with general monotone nonlinear term and possibly
degenerate noise. Non-degeneracy is replaced by a structure
condition as it usually happens in BSDEs approach, see, for
instance, \cite{ElKaMaz}, \cite{FuTe1}. Moreover the use of
$L^{\infty}$ estimates specific to infinite horizon backward
stochastic differential equations (see \cite{bh}, \cite{royer},
\cite{HuTe}) allow us to eliminate conditions on  the lipschitz
constant of the hamiltonian. On the other side we will only
consider bounded cost functionals.

$ $

 To start being more precise we consider a forward equation
$$dX_t^x=(AX_t^x+F(X_t^x))dt+G dW_t,\qquad X_0=x$$
where $X$ has values in a Banach space $E$, $F$ maps $E$ to $E $
and $A$ generates a strongly continuous semigroup of contractions.
Appropriate dissipativity assumptions on  $A+F$   ensure the
exponential decay of the difference between the trajectories
starting from different points $x,x'\in E$.

Then we introduce the class of strictly monotonic backward
stochastic differential equations
\begin{equation}\label{bsderoyer*}
{Y}^{x,\alpha}_t={Y}^{x,\alpha}_T +\int_t^T(\psi(X^{x}_\sigma,
Z^{x,\alpha}_\sigma)-\alpha Y^{x,\alpha}_\sigma)d\sigma-\int_t^T
Z^{x,\alpha}_\sigma dW_\sigma, \quad 0\le t\le T <\infty.
\end{equation}
for all $\alpha>0$ (see \cite{bh}, \cite{royer} or \cite{HuTe})
where $\psi: E\times\Xi^*\rightarrow \mathbb{R}$ is bounded in the
first variable and Lipschitz in the second. By estimates based on
a Girsanov argument introduced in \cite{bh} we obtain uniform
estimates on $\alpha{Y}^{x,\alpha}$ and
${Y}^{x,\alpha}-{Y}^{x',\alpha}$ that allow us to prove that,
roughly speaking, $({Y}^{x,\alpha}-{Y}^{0,\alpha}_0,
{Z}^{x,\alpha}, \alpha {Y}^{0,\alpha}_0)$ converge to a solution
$(Y^x,Z^x,\lambda)$ of the EBSDE (\ref{EBSDE*}), for all $x\in E$.
We also show that $\lambda$ is unique under very general
conditions. On the contrary, in general we can not  expect
uniqueness of the solution to (\ref{EBSDE*}), at least in the non
markovian case. On the other side in the markovian case we show
that we can find a solution of (\ref{EBSDE*}) with
$Y^x_t=v(X^x_t)$ and $Z^x_t=\zeta(X^x_t)$ where $v$ is Lipschitz
and $v(0)=0$. Moreover $(v, \zeta)$ are unique at least in a
special case where $\psi$ is the Hamiltonian of a control problem
and the processes $X^x$ are recurrent (see Section \ref{sec-uniq}
where we adapt an argument from \cite{GoMa}).

$ $

If we further assume differentiability  of $F$ and
 $\psi$ (in the Gateaux sense) then $v$ is differentiable,
  moreover $\zeta =\nabla v G$
 and finally $(v,\lambda)$ give a mild solution of the HJB equation
 \begin{equation}
\mathcal{L}v(x)
+\psi\left( x,\nabla v(x)  G\right)  = \lambda, \quad x\in E,  \label{hjb*}%
\end{equation}
where  linear operator $\mathcal{L}$ is formally defined by
\[
\mathcal{L}f\left(  x\right)  =\frac{1}{2}Trace\left(
GG^{\ast}\nabla ^{2}f\left(  x\right)  \right)  +\langle Ax,\nabla
f\left(  x\right) \rangle_{E,E^{\ast}}+\langle F\left(  x\right)
,\nabla f\left(  x\right) \rangle_{E,E^{\ast}}.
\]
Moreover if the Kolmogorov semigroup satisfies the smoothing
property in Definition \ref{strongly-feller} and $F$ is genuinely
dissipative (see Definition \ref{gen-diss}) then $v$ is bounded.

$ $

The above results are then applied to a control problem with cost
\begin{equation}\label{ergodic-cost*}
J(x,u)=\limsup_{T\rightarrow\infty}\frac{1}{T}\, \mathbb E\int_0^T
L(X_s^x,u_s)ds,
\end{equation}
 where $u$ is an adapted process (an admissible control)
with values in a separable metric space $U$, and the state
equation is a Banach valued evolution equation of the form
$$dX_t^x=(AX_t^x+F(X_t^x))\, dt+G(dW_t+R(u_t)\,dt),$$
where $R: U \rightarrow \Xi$ is bounded. It is clear that the
above functional depends only on the asymptotic behavior of the
trajectories of $X^x$. After appropriate formulation
 we prove that, setting $\psi(x,z)= \inf_{u\in U} [L(x,u)+ zR(u)]$ in
  (\ref{EBSDE*}), then $\lambda$ is optimal, that is
  $$\lambda=\inf_{u}J(x,u)$$
   where the infimum is over all admissible controls.
   Moreover $Z$ allows to  construct on optimal feedback in the
   sense that $$\lambda=J(x,u) \hbox{ if and only if } L(X_t^x,u_t)+Z_t
R(u_t)=\psi(X_t^x,Z_t).$$

Finally, see Section \ref{section-heat-eq}, we show that our
assumptions allow us to treat ergodic optimal control problems for
a stochastic heat equation with polynomial nonlinearity and
space-time white noise. We notice that the Banach space setting is
essential in order to treat nonlinear terms with superlinear
growth  in the state equation.

$ $

The paper is organized as follows.
After a section on notation, we introduce the forward SDE; in section 4 we
study the ergodic BSDEs; in section 5 we show in
addition the differentiability of the
solution assuming that the coefficient is Gateaux differentiable.
In section 6 we study the ergodic Hamilton-Jacobi-Bellman
equation and we apply our result to optimal
ergodic control in section 7. Section 8 is devoted to
show the uniqueness of Markovian solution and the last section
contains application to the ergodic control of a nonlinear stochastic
heat equation.

\section{Notation}
 Let $E,F$ be Banach spaces, $H$ a Hilbert space, all
assumed to be defined over the real field and to be separable. The
norms and the scalar product will be  denoted $|\,\cdot\,|$,
$\langle\,\cdot\,,\,\cdot\,\rangle$, with subscripts if needed.
Duality between the dual space $E^*$ and
 $E$ is denoted $\langle\,\cdot\,,\,\cdot\,\rangle_{E^*,E}$.
$L(E,F)$ is the space of linear bounded operators $E\to F$, with
the operator norm.
The domain of a linear (unbounded) operator $A$ is denoted $D(A)$.

Given a bounded function
$ \phi: E\rightarrow \mathbb{R}$ we denote
$\Vert\phi\Vert_0=\sup_{x\in E}|\phi(x)|$. If, in addition,
$\phi$ is also Lipschitz continuous then
$\Vert\phi\Vert_{\hbox{lip}}=\Vert\phi\Vert_0+
\sup_{x,x'\in E,\,x\ne x'}|\phi(x)-\phi(x')||x-x'|^{-1}$.

We say that a function $F:E\to F$ belongs to
the class $\calg^1(E,F)$ if it is continuous, has a Gateaux
differential $\nabla F(x)\in L(E,F)$ at any point $x  \in E$, and
for every $k\in E$ the mapping $x\to \nabla F(x) k$  is continuous
from $E$ to $F$ (i.e. $x\to \nabla F(x) $  is continuous from $E$
to $L(E,F)$ if the latter space is endowed the strong operator
topology). In connection with stochastic equations,
the space $\calg^1$ has been introduced in \cite{FuTe1},
to which we refer the reader for further properties.

 Given a  probability space $\left(
\Omega,\mathcal{F},\mathbb{P}\right) $ with a filtration
$(\calf_t)_{t\ge 0}$ we consider the following classes of
stochastic processes with values in a real separable Banach space
$K$.

\begin{enumerate}
\item
$L^p_{\mathcal{P}}(\Omega,C([0,T],K))$, $p\in [1,\infty)$,
$T>0$, is the space
of predictable processes $Y$ with continuous paths
on $[0,T]$
such that
$$
|Y|_{L^p_{\mathcal{P}}(\Omega,C([0,T],E))}^p
= \E\, \sup_{t\in [0,T]}|Y_t|_K^p<\infty.
$$
\item
$L^p_{\mathcal{P}}(\Omega,L^2([0,T];K))$, $p\in [1,\infty)$,
$T>0$, is the space
of predictable processes $Y$ on $[0,T]$ such that
$$
|Y|^p_{L^p_{\mathcal{P}}(\Omega,L^2([0,T];K))}=
\E\,\left( \int_{0}^{T}|Y_t|_K^2\,dt\right)^{p/2}<\infty.
$$

\item
$L_{\cal P, {\rm loc}}^2(\Omega;L^2(0,\infty;K))$
is the space
of predictable processes $Y$ on $[0,\infty)$ that belong
to the space $L^2_{\mathcal{P}}(\Omega,L^2([0,T];K))$
for every $T>0$.
\end{enumerate}

\section{The forward equation}
In a complete probability space $\left(
\Omega,\mathcal{F},\mathbb{P}\right) ,$  we consider the following
stochastic differential equation with values in a  Banach
space $E$:
\begin{equation}
\left\{
\begin{array}[c]{l} dX_t  =AX_t  d t+F(X_t) dt +GdW_t ,\text{ \ \ \
} t\geq 0, \\
X_0  =x\in \, E.
\end{array}
\right.  \label{sde}
\end{equation}
We assume that $E$ is continuously and densely embedded in a
Hilbert space $H$, and that both spaces are real separable.

 We will work under the following
general assumptions:

\begin{hypothesis}
\label{general_hyp_forward}
\begin{enumerate}
 \item The operator $A$ is the generator of  a strongly
 continuous semigroup of contractions in
$E$. We assume that the semigroup
 $\{e^{tA},\, t\geq0\}$
 of bounded linear operators on $E$ generated by $A$
 admits an extension to a strongly
continuous semigroup of bounded linear operators on $H$ that we
denote by $\{S(t),\, t\geq 0\}$.

\item $W$ is a cylindrical Wiener process in another real separable
Hilbert space $\Xi$. Moreover by $\calf_{t}$ we denote the
$\sigma$-algebra generated by $\{W_s,\; s\in [0,t]\}$ and
 by the sets of $\calf$ with $\P$-measure zero.

\item $F:E\to E$ is continuous and has polynomial growth (that is
there exist $c>0, k\ge 0$ such that $|F(x)|\leq c (1+|x|^k)$,
$x\in E$). Moreover there exists $\eta>0$
 such that $A+F+\eta I$ is dissipative.

\item $G$ is a bounded linear operator from $\Xi$ to $H$. The
bounded linear, positive and symmetric  operators on $H$ defined
by the formula
\[
Q_{t}h=\int_{0}^{t}S(s)GG^{\ast}S^*(s)h\,ds,\qquad t\geq 0,\; h\in
H,
\]
are assumed to be of trace class in $H$. Consequently we can
define the stochastic convolution
$$
W^{A}_t  =\int_{0}^{t}S(t-s) GdW_s,\quad t\geq 0,
$$
as a family of $H$-valued stochastic integrals. We assume that the
process $\{W^{A}_t,\, t \geq 0\}$ admits an $E$-continuous
version.
\end{enumerate}
\end{hypothesis}

We recall that, for every $x\in E$, with $x\neq 0$, the
subdifferential of the norm at $x$, $\partial\left( |x| \right) $,
is the set of functionals $x^{\ast}\in E^{\ast}$ such that
$\left\langle x^{\ast },x\right\rangle _{E^{\ast},E}=| x| $ and $|
x^{\ast}|_{E^{\ast}}=1$. If $x=0$ then $\partial\left( | x|\right)
$ is the set of functionals $x^{\ast}\in E^{\ast}$ such that
$|x^{\ast}|_{E^{\ast}}\leq 1$. The dissipativity assumption on
$A+F$ can be explicitly stated as follows: for $x,x'\in
D(A)\subset E$ there exists $x^{\ast } \in\partial\left( \left|
x-x' \right| \right)  $ such that
$$\left\langle
x^{\ast}  ,A( x-x'  )  +F\left( x \right) -F\left( x' \right)
\right\rangle _{E^{\ast},E}
 \leq-\eta\left|
x-x' \right|.
$$

We can state the following theorem, see e.g. \cite{DP1}, theorem
7.13 and \cite{DP2}, theorem 5.5.13.
\begin{theorem}
\label{teo2 forward}Assume that Hypothesis
\ref{general_hyp_forward} holds true. Then for every $x\in E$
equation (\ref{sde}) admits a unique mild solution, that is an
adapted  $E$-valued process with continuous paths  satisfying
$\mathbb{P}$-a.s.
\[
X_{t}=e^{  t  A}x+\int_{0}^{t}e^{ (t-s ) A}F\left(  X_{s}\right)
ds+\int_{0}^{t}e^{(t-s ) A}GdW_{s},\text{ \ \ \ }t\geq 0 .
\]
\end{theorem}

We denote the solution by $X^x  $, $x\in E$.

 Now we want to investigate the dependence of
the solution on the initial datum.

\begin{proposition}
\label{prop lip X} Under Hypothesis \ref{general_hyp_forward}  it
holds:
\[
\left| X_t^{x_1}  -X_t^{x_2} \right| \leq e^{-\eta t  }\left|
x_{1}-x_{2}\right| ,\text{ }t\ge 0, \;\; x_{1},x_{2}\in E.
\]

\end{proposition}

\begin{proof}
Let $X_{1}\left(  t\right)  =X^{x_1}_{t} $ and $X_{2}\left(
t\right)  =X^{x_2}_{t} $, $x_{1},x_{2}\in E$. For $i=1,2$ we set
$X_{i}^{n}\left(  t\right) =J_n X_{i}\left(  t\right)  $, where
$J_n
 =n\left(  nI-A\right) ^{-1}$. Since $X_{i}^{n}\left(  t\right)
\in {D}\left(  A\right)  $ for every $t\geq 0 $, and
\[
X_{i}^{n}\left(  t\right)  =e^{t A}J_n x_{i}+\int
_{0}^{t}e^{\left(  t-s\right)  A}J_n F\left( X_{i}\left(  s\right)
\right) ds+\int_{0}^{t}e^{\left(  t-s\right) A}J_n GdW_{s},
\]
we get
\[
\frac{d}{dt}\left(  X_{1}^{n}\left(  t\right) -X_{2}^{n}\left(
t\right)  \right)  =A\left(  X_{1}^{n}\left(  t\right)  -X_{2}
^{n}\left(  t\right)  \right)  +J_n \left[ F\left( X_{1}\left(
t\right)  \right)  -F\left( X_{2}\left( t\right)  \right)  \right]
.
\]
So, by proposition II.8.5 in \cite{S}\ also $\left|
X_{1}^{n}\left( t\right)  -X_{2}^{n}\left(  t\right) \right| $
admits the left and right derivatives with respect to $t$ and
there exists $x_{n}^{\ast }\left(  t\right) \in\partial\left(
\left| X_{1}^{n}\left( t\right) -X_{2}^{n}\left(  t\right) \right|
\right)  $ such that the left derivative of $\left|
X_{1}^{n}\left( t\right) -X_{2}^{n}\left(  t\right)  \right| $
satisfies the following
\[
\frac{d^{-}}{dt}\left| X_{1}^{n}\left(  t\right) -X_{2}^{n}\left(
t\right)  \right| =\left\langle x_{n}^{\ast}\left(  t\right)
,\frac{d}{dt}\left( X_{1}^{n}\left(  t\right)  -X_{2}^{n}\left(
t\right)  \right) \right\rangle _{E^{\ast},E}.
\]
So we have
$$ \begin{array}{ll}\dis \frac{d^{-}}{dt}\left| X_{1}^{n}\left(
t\right) -X_{2}^{n}\left( t\right)  \right|   & =\left\langle
x_{n}^{\ast}\left(  t\right)  ,A\left( X_{1}^{n}\left( t\right)
-X_{2}^{n}\left(  t\right)  \right)  +F\left(  X_{1} ^{n}\left(
t\right)  \right)  -F\left(  X_{2}^{n}\left( t\right)
\right)  \right\rangle _{E^{\ast},E}\\
& \quad +\left\langle x_{n}^{\ast}\left(  t\right)  ,J_n F\left(
X_{1}\left(  t\right)  \right)  -F\left(  X_{1}
^{n}\left(  t\right)  \right)  \right\rangle _{E^{\ast},E}\\
&  \quad -\left\langle x_{n}^{\ast}\left(  t\right)  ,J_n F\left(
X_{2}\left(  t\right)  \right)  -F\left(  X_{2}
^{n}\left(  t\right)  \right)  \right\rangle _{E^{\ast},E}\\
&  \leq-\eta\left|
 X_{1}^{n}\left(  t\right)  -X_{2}^{n}\left(
t\right)  \right| +\left| \delta_{1}^{n}\left( t\right)
-\delta_{2}^{n}\left(  t\right)  \right| ,
\end{array}$$
where for $i=1,2$ we have set $\delta_{i}^{n}\left( t\right) =J_n
F\left(  X_{i}\left( t\right) \right) -F\left( X_{i}^{n}\left(
t\right)  \right)  $.

Multiplying the above by $e^{\eta t}$ we get
$$\frac{d^{-}}{dt}\left( e^{\eta t}\left| X_{1}^{n}\left( t\right)
-X_{2}^{n}\left( t\right) \right|\right)\leq e^{\eta t} \left|
\delta_{1}^{n}\left( t\right) -\delta_{2}^{n}\left( t\right)
\right|.$$ We note that $\delta_{i} ^{n}\left(  t\right)  $ tends
to $0$\ uniformly in $t\in \left[0,T\right]  $ for arbitrary
$T>0$. Indeed,
\[
\delta_{i}^{n}\left(  t\right)  =nR\left(  n,A\right)  \left[
F\left( X_{i}\left(  t\right)  \right)  -F\left( X_{i}^{n}\left(
t\right)  \right)  \right]  +\left( nR\left( n,A\right) -I\right)
F\left(  X_{i}\left( t\right) \right)  ,
\]
and the convergence to $0$ follows by a classical argument, see
e.g. the proof of theorem 7.10 in \cite{DP1}, since
$X_{i}^{n}\left( t\right)  $ tends to $X_{i}\left(  t\right)  $
uniformly in $t\in\left[  0,T\right]  $ and the maps $t\mapsto
X_{i}\left(  t\right)   $ and $t\mapsto F\left( X_{i}\left(
t\right)  \right)  $ are continuous with respect to $t$.

 Thus letting  $n\rightarrow\infty$ we can conclude
\[
\left| X_{1}\left(  t\right)  -X_{2}\left(  t\right) \right| \leq
e^{-\eta t } \left| x_{1}-x_{2}\right| .
\]
and the claim is proved. \end{proof}

$ $

\noindent We will also need the following assumptions.
\begin{hypothesis}
\label{hyp_W_A F(W_A)} We have $\sup_{t\geq 0}\,
\E\,|W^A_t|^2<\infty.$
\end{hypothesis}
\begin{hypothesis}
\label{hyp-convol-determ}  $e^{tA}G\,(\Xi)\subset E$ for all $t>0$
and $\displaystyle \int_0^{+\infty} |e^{tA} G|_{L(\Xi,E)} dt <
\infty$.
\end{hypothesis}

We recall that for arbitrary gaussian random variabile $Y$ with
values in the Banach space $E$, the inequality
$$
\E \,\phi (|Y|-\E\,|Y|)\le \E \,\phi (2\sqrt{\E\,|Y|^2}\,\gamma)
$$
holds for any convex nonnegative continuous function $\phi$
on $E$ and
for $\gamma$ a real standard gaussian random variable, see e.g.
\cite{kw-woy}, Example 3.1.2. Upon taking $\phi(x)=|x|^p$,  it
follows that for every $p\ge 2$ there exists $c_p>0$ such that $\E
\,|Y|^p\le c_p(\E \,|Y|^2)^{p/2}$. By the gaussian character of
$W^A_t$ and the polynomial growth condition on $F$ stated in
Hypothesis \ref{general_hyp_forward}, point 3, we see that
Hypothesis \ref{hyp_W_A F(W_A)} entails that for every $p\ge 2$
\begin{equation}\label{stimegaussunif}
\sup_{t\geq 0} \E\left[ |W^A_t|^p+ |F(W^A_t)|^p \right] <\infty.
\end{equation}

\begin{proposition}\label{prop-X-L^p}
Under Hypothesis \ref{general_hyp_forward} it holds, for arbitrary
$T>0$ and arbitrary $p\geq 1$
\begin{equation}\label{prop-X-L^p-1}
\E\sup_{t\in [0,T]}  |X_t^x|^p \leq C_{p,T}(1+|x|^p),\qquad x\in
E.
\end{equation}
 If, in addition, Hypothesis
\ref{hyp_W_A F(W_A)} holds then, for a suitable constant C
\begin{equation}\label{prop-X-L^p-2}\sup_{t\geq 0} \E |X_t^x| \leq C(1+|x|)
,\qquad x\in E.\end{equation} Moreover if, in addition, Hypothesis
\ref{hyp-convol-determ} holds, $\gamma$ is a bounded, adapted,
$\Xi$-valued process and $X^{x,\gamma}$ is the mild solution of
equation
\begin{equation}
\left\{
\begin{array}{l}
dX^{x,\gamma}_t  =AX^{x,\gamma}_{t} dt+F(
X^{x,\gamma}_{t} ) dt+GdW_{t}+G\gamma_{t}\,dt ,\quad t\geq 0, \\
X^{x,\gamma}_{0} =x\in E.
\end{array} \right.  \label{sde-gamma}
\end{equation}
then it is still true that
\begin{equation}\label{rel-estimate-Xgamma}
 \sup_{t\geq 0} \E |X^{x,\gamma}_t| \leq C_{\gamma}(1+|x|),\qquad  x\in E,
\end{equation}
for a suitable constant $C_{\gamma}$ depending only on
 a uniform bound for $\gamma$.
\end{proposition}
\begin{proof} We let $Z_t=X^x_t-W^A_t$,
$Z^n_t=J_n Z_t $, then
$$\frac{d}{dt } Z^n_t =
AZ^n_t +J_nF(X^x_t) = AZ^n_t +\left[F(Z^n_t+J_n W^A_t) - F(J_n
W^A_t)\right]+F( W^A_t)+\delta^n_t
$$ where $$\delta^n_t= J_n F(X^x_t)-F(J_n X^x_t)
+F(J_n W^A_t)-F( W^A_t).$$ Proceeding as in the proof of
Proposition  \ref{prop lip X} observing  that, for all $t>0$,
 $\displaystyle \int_0^{t}|\delta^n_s| ds \rightarrow 0$ as $n\rightarrow\infty$, we get:
$$|Z_t|\leq e^{-\eta t}|x|+\int_0^{t} e^{-\eta (t-s)}
|F(W^A_s)|ds,\;\;\; \mathbb{P}-\hbox{a.s.}$$ and
(\ref{prop-X-L^p-2}) follows from (\ref{stimegaussunif}).

In the case in which $X^x$ is replaced by $X^{x, \gamma}$ the
proof is exactly the same just replacing $W^A_t$ by
$W^{A,\gamma}_t=W^A_t+\int_0^t e^{(t-s)A}G\gamma_s ds$.

Finally to prove (\ref{prop-X-L^p-1}) we notice that (see the
discussion in \cite{masiero}) the process $W^A$ is a Gaussian
random variable with values in $C([0,T],E)$. Therefore by the
polynomial growth of $F$ we get
$$ \E\sup_{t\in [0,T]} \left[|W^A_t|^p + |F(W^A_t)|^p\right]\leq
C_{p,T}(1+|x|^p),$$ and the claim follows as above.
\end{proof}

$ $

 Finally the following result  is proved exactly as
Theorem 6.3.3. in \cite{DP2}.
\begin{theorem}\label{ergodicity}
Assume that Hypotheses \ref{general_hyp_forward} and
 \ref{hyp_W_A F(W_A)} hold then equation (\ref{sde}) has a unique
 invariant measure in $E$ that we will denote by $\mu$. Moreover
 $\mu$ is strongly
mixing (that is, for all $x\in E$, the law of $X_t^x$ converges
weakly to $ \mu$ as
 $t\rightarrow \infty$).
 Finally
there exists a constant $C>0$ such that for any bounded Lipschitz
function $\phi: E\rightarrow \mathbb{R}$,
$$\left|\mathbb{E}\phi(X^x_t)-\int_E \phi\, d\mu \right|\leq C(1+|x|)
e^{-\eta t /2} \Vert\phi\Vert_{\hbox{\em lip}}.$$
\end{theorem}

\section{Ergodic BSDEs (EBSDEs)}

This section is devoted to the following type of BSDEs with
infinite horizon
\begin{equation}\label{EBSDE}
Y^x_t=Y^x_T +\int_t^T\left[\psi(X^x_\sigma,Z^x_\sigma)-
\lambda\right]d\sigma-\int_t^T Z^x_\sigma\, dW_\sigma, \quad 0\le
t\le T <\infty,
\end{equation}
where $\lambda$ is a real number and is part of the unknowns of
the problem; the equation is required to hold for every $t$ and
$T$ as indicated. On the function $\psi: E\times  \Xi^*
\rightarrow {\mathbb R}$ and assume the following:

\begin{hypothesis}\label{hypothesisroyer} $ $ There exists
$K_x, K_z>0$ such that
$$ |\psi(x,z)
-\psi(x',z')|\le K_x|x-x'|+ K_z |z-z'|, \qquad
 x,x'\in E,\;
z,z'\in\Xi^*.
$$
Moreover $\psi(\,\cdot\,,0)$ is bounded.  We denote $\sup_{x\in E
}|\psi(x,0)|$ by $M$.
\end{hypothesis}
We start by considering an infinite horizon equation with strictly
monotonic drift, namely, for $\alpha>0$, the equation
\begin{equation}\label{bsderoyer}
{Y}^{x,\alpha}_t={Y}^{x,\alpha}_T +\int_t^T(\psi(X^{x}_\sigma,
Z^{x,\alpha}_\sigma)-\alpha Y^{x,\alpha}_\sigma)d\sigma-\int_t^T
Z^{x,\alpha}_\sigma dW_\sigma, \quad 0\le t\le T <\infty.
\end{equation}

The existence and uniqueness of solution to (\ref{bsderoyer})
under Hypothesis \ref{hypothesisroyer} was first studied by Briand
and Hu in \cite{bh} and then generalized by Royer in \cite{royer}.
 They have established the following result when $W$ is a finite dimensional
  Wiener process but the extension to the case in which $W$ is a
  Hilbert-valued Wiener process is immediate (see also  \cite{HuTe}).

\begin{lemma}\label{lemmaroyer} Let us suppose that Hypotheses
\ref{general_hyp_forward} and
\ref{hypothesisroyer} hold.
 Then
 there  exists a unique solution $(Y^{x,\alpha},Z^{x,\alpha})$
 to  BSDE (\ref{bsderoyer})
such that $Y^{x,\alpha}$ is a bounded continuous process, and
$Z^{x,\alpha}$ belongs to $L_{\cal P, {\rm
loc}}^2(\Omega;L^2(0,\infty;\Xi^*))$.

Moreover $|Y^{x,\alpha}_t|\leq {M}/{\alpha}$, $\mathbb{P}$-a.s.
for all $t\geq 0$.
\end{lemma}
We define $$v^{\alpha}(x)=Y^{\alpha,x}_0.
$$
We notice that by the above $|v^{\alpha}(x)|\leq {M}/{\alpha}$ for
all $x\in E$. Moreover by the uniqueness of the solution of
equation (\ref{bsderoyer}) it follows that
$Y^{\alpha,x}_t=v^{\alpha}(X^x_t)$

To establish Lipschitz continuity of $  v^{\alpha}$ (uniformly in
$\alpha$)  we use a Girsanov argument due to P. Briand and Y. Hu,
see \cite{bh}. Here and in the following we use an
infinite-dimensional version of the Girsanov formula that can be
found for instance in \cite{DP1}.
\begin{lemma}\label{lemma-lip-v}  Under Hypotheses \ref{general_hyp_forward}
and \ref{hypothesisroyer}  the following
holds for any $\alpha>0$:
$$|v^{\alpha}(x) - v^{\alpha}(x')| \leq \frac{K_x}{\eta} |x-x'|,
\qquad  x,x'\in E. $$
\end{lemma}
\begin{proof} We briefly report the argument for the reader's convenience.

We set $\tilde{Y}=Y^{\alpha,x}-Y^{\alpha,x'}$,
$\tilde{Z}=Z^{\alpha,x}-Z^{\alpha,x'},$
$$\beta_t=\begin{cases}
\frac{\displaystyle \psi(X^{x'}_t,Z^{\alpha,
x'}_t)-\psi(X^{x'}_t,Z^{\alpha,x}_t)} {\displaystyle
|Z^{\alpha,x}_t - Z^{\alpha,x'}_t|_{\Xi^*}^2}\left( Z^{\alpha,x}_t
- Z^{\alpha,x'}_t\right)^*,& \hbox{ if } Z^{\alpha,x}_t \neq Z^{\alpha,x'}_t \\
0, & \hbox{  elsewhere,   }
 \end{cases}
$$
$$f_t=\psi(X^{x}_t,
Z^{x,\alpha}_t)-\psi(X^{x'}_t, Z^{x,\alpha}_t).    $$ By
Hypothesis  \ref{hypothesisroyer}, $\beta$ is a bounded
$\Xi$-valued, adapted process thus there exists a probability
$\tilde{\mathbb{P}}$ under which $\tilde{W_{t}}=\int_0^{t} \beta_s
ds + W_{t}$ is a cylindrical $\Xi$-valued Wiener process for
${t}\in [0,T]$. Then  $(\tilde{Y},\tilde{Z})$ verify, for all
$0\le t\le T <\infty$,
\begin{equation}\label{bsderoyer-girsanov}
\tilde{Y}_t=\tilde{Y}_T -\alpha \int_t^T \tilde{Y}_\sigma d\sigma
+\int_t^T f_{\sigma}d\sigma- \int_t^T \tilde{Z}_\sigma
d\tilde{W}_\sigma.
\end{equation}
Computing  $d (e^{-\alpha t}\tilde{Y}_t)$, integrating over
$[0,T]$, estimating the absolute value and finally taking the
conditional expectation
 $\tilde{\mathbb{E}}^{\mathcal{F}_t}$ with respect to
$\tilde{\mathbb{P}}$ and $\mathcal{F}_t$ we get:
$$  |\tilde{Y}_t| \leq  e^{-\alpha(T-t)}  \tilde{\mathbb{E}}^{\mathcal{F}_t}
|   \tilde{Y}_T |+
 \tilde{\mathbb{E}}^{\mathcal{F}_t}
  \int_{t}^T  e^{-\alpha(s-t)} |f_s| ds  $$
Now we recall that $ \tilde{Y}$ is bounded   and that $|f_t|\leq
K_x |X^{x}_t-X^{x'}_t|\leq K_x e^{-\eta t}|x-x'|$ by Proposition
\ref{prop lip X}. Thus if  $T\rightarrow \infty$ we get $
|\tilde{Y}_t| \leq K_x (\eta+\alpha)^{-1}e^{\alpha t} |x-x'|  $
and the claim follows  setting $t=0$.
\end{proof}

$ $

\noindent By the above Lemma if we set
$$\overline{v}^{\alpha}(x)= {v}^{\alpha}(x)- {v}^{\alpha}(0),$$
then $   | \overline{v}^{\alpha}(x)|\leq    K_x \eta^{-1}|x|$ for
all $x\in E$ and all $\alpha>0$. Moreover by Lemma
\ref{lemmaroyer} $\alpha |{v}^{\alpha}(0)|\leq M$.

 \noindent  Thus by a diagonal procedure   we can construct a
 sequence $\alpha_n\searrow 0$ such that for all $x$ in a
countable dense subset $D\subset E$
     \begin{equation}\label{def-of-lambda}
    {\overline{v}}^{\alpha_n}(x)\rightarrow \overline{v}(x),\qquad
\alpha_n v^{\alpha_n}(0)\rightarrow \overline{\lambda},
    \end{equation}
for a suitable  function $  \overline{v}: D \rightarrow
\mathbb{R}$ and for a suitable real number $\overline{\lambda}$.

 Moreover, by Lemma \ref{lemma-lip-v}, $ | \overline{v}^{\alpha}(x)-
\overline{v}^{\alpha}(x')|\leq    K_x \eta^{-1}|x-x'|$   for all
$x,x'\in E$ and all $\alpha>0$. So $\overline{v}$ can be extended
to a Lipschitz function defined on the whole $E$ (with Lipschitz
constant $K_x \eta^{-1}  $) and
\begin{equation}\label{def-of-v} {\overline{v}}^{\alpha_n}(x)\rightarrow
\overline{v}(x),\qquad  x\in E.\end{equation}

\begin{theorem} \label{main-EBSDE} Assume Hypotheses
\ref{general_hyp_forward}  and
\ref{hypothesisroyer} hold. Moreover let $\bar \lambda$ be the
real number in (\ref{def-of-lambda}) and define $\bar Y^x_t= \bar
v(X^x_t)$ (where $\overline{v}$ is the Lipschitz function with
$\overline{v}(0)=0$ defined in (\ref{def-of-v})). Then there
exists a process $\overline{Z}^{x}\in L_{\cal P, {\rm
loc}}^2(\Omega;L^2(0,\infty;\Xi^*))$
  such that $\mathbb{P}$-a.s.  the EBSDE
 (\ref{EBSDE}) is satisfied by
 $(\bar Y^x,\bar Z^x, \bar \lambda)$ for all $0\leq t\leq T$.

Moreover  there exists a measurable function  $\overline{\zeta}:
E\rightarrow \Xi^*$ such that
$\overline{Z}^{x}_t=\overline{\zeta}(X^x_t)$.
\end{theorem}

\begin{proof} Let $\overline{Y}^{x,\alpha}_t={Y}^{x,\alpha}_t-v^{\alpha}(0)=
\overline{v}^{\alpha}({X}^{x}_t)$.  Clearly we have,
$\mathbb{P}$-a.s.,
\begin{equation}\label{equation-proof-main-1}
 \overline{Y}^{x,\alpha}_t=\overline{Y}^{x,\alpha}_T +\int_t^T(\psi(X^{x}_\sigma,
Z^{x,\alpha}_\sigma)-\alpha \overline{Y}^{x,\alpha}_\sigma-\alpha
{v}^{\alpha}(0))d\sigma -\int_t^T Z^{x,\alpha}_\sigma dW_\sigma,
\quad 0\le t\le T <\infty.
\end{equation}
Since $|\bar v^{\alpha}(x)|\leq K_x|x|/\eta $, inequality
(\ref{prop-X-L^p-1})   ensures that
$\mathbb{E}\sup_{t\in[0,T]}\left[\sup_{\alpha>0}
|\overline{Y}^{x,\alpha}_t|^2\right]< +\infty$  for any $T>0$.
Thus, if we define $\overline{Y}^x=\overline{v}(X^x)$, then by
dominated convergence theorem
$$\mathbb{E} \int_0^T |\overline{Y}^{x,\alpha_n}_t -\overline{Y}^{x}_t|^2 dt
 \rightarrow 0\quad \hbox{and}\quad
\mathbb{E} |\overline{Y}^{x,\alpha_n}_T-\overline{Y}^{x}_T|^2
\rightarrow 0
$$
as $n\rightarrow \infty$ (where $\alpha_n \searrow 0$ is a
sequence for which (\ref{def-of-lambda}) and (\ref{def-of-v})
hold).

We claim now that there exists $\overline{Z}^{x}\in L_{\cal P,
{\rm loc}}^2(\Omega;L^2(0,\infty;\Xi^*))$ such that
 $$\mathbb{E} \int_0^T |{Z}^{x,\alpha_n}_t -\overline{Z}^{x}_t|_{\Xi^*}^2 dt
  \rightarrow 0$$
Let  $\tilde{Y}={\bar Y}^{x,\alpha_n}-{\bar Y}^{x,\alpha_m}$,
$\tilde{Z}={Z}^{x,\alpha_n}-{Z}^{x,\alpha_m}$. Applying It\^o's
rule to $\tilde{Y}^2$ we get by standard computations
$$\tilde{Y}^2_0+\mathbb{E}\int_0^T |\tilde{Z}_t|_{\Xi^*}^2 dt
=\mathbb{E}{\tilde Y}^2_T + 2\mathbb{E}\int_0^T \tilde \psi_t
\tilde Y_t dt -2 \mathbb{E}\int_0^T \left[\alpha_n
{Y}^{x,\alpha_n}_t - \alpha_m {Y}^{x,\alpha_m}_t\right] \tilde
Y_t\,dt
$$
where  $\tilde
\psi_t=\psi(X^x_t,Z^{x,\alpha_n}_t)-\psi(X^x_t,Z^{x,\alpha_m}_t)$.
We notice that $|\tilde\psi_t| \leq  K_z|\tilde Z _t|$ and
$\alpha_n |{Y}^{x,\alpha_n}_t|\leq M$. Thus
$$
\mathbb{E}\int_0^T |\tilde{Z}_t|_{\Xi^*}^2 dt \leq c\left[
\mathbb{E} (\tilde Y^x_T)^2 +\mathbb{E}\int_0^T (\tilde{Y}^x_t)^2
dt +\mathbb{E}\int_0^T |\tilde{Y}^x_t| dt \right].$$ It follows
that the sequence $\{{Z}^{x,\alpha_m}\}$ is Cauchy in
$L^2(\Omega;L^2(0,T;\Xi^*))$ for all $T>0$ and our claim is
proved.

Now we can  pass to the limit as $n\rightarrow \infty$ in equation
(\ref{equation-proof-main-1}) to obtain
\begin{equation}\label{equation-proof-main-2}
 \overline{Y}^{x}_t=\overline{Y}^{x}_T +\int_t^T(\psi(X^{x}_\sigma,
\overline{Z}^{x}_\sigma)-\overline{\lambda })d\sigma-\int_t^T
\overline{Z}^{x}_\sigma dW_\sigma, \quad 0\le t\le T <\infty.
\end{equation}
We notice  that the above equation also ensures continuity of the
trajectories of $\overline{Y}$ It remains now to prove that we can
find a measurable function $\bar \zeta:E\rightarrow \Xi^*$ such
that
 $\overline{Z}^{x}_t=\bar \zeta (X^x_t)$, $\mathbb{P}$-a.s. for almost every $t\geq 0$.

By a general argument, see for instance \cite{Fu}, we know that
for all $\alpha>0$ there exists  $\zeta^{\alpha}:E\rightarrow
\Xi^*$ such that
 ${Z}^{x,\alpha}_t=\zeta^{\alpha} (X^x_t)$, $\mathbb{P}$-a.s.
 for almost every $t\geq 0$.

To construct $\zeta$ we need some more regularity of the processes
${Z}^{x,\alpha}$ with respect to $x$.

If we compute $d ({Y}^{x,\alpha}_t-{Y}^{x',\alpha}_t)^2$ we get by
the Lipschitz character of $\psi$:
$$ \begin{array} {l}
\displaystyle \mathbb{E}\int_0^T
|Z^{x,\alpha}_t-Z^{x',\alpha}_t|_{\Xi^*}^2 dt \leq \mathbb{E}
(v^{\alpha}(X^x_T)- v^{\alpha}(X^{x'}_T))^2
  \\
\quad +  \displaystyle \mathbb{E}\int_0^T
\left(K_x|X^x_s-X^{x'}_s|
+K_z|Z^{x,\alpha}_s-Z^{x',\alpha}_s|\right)
\left|v^{\alpha}(X^x_s)- v^{\alpha}(X^{x'}_s)\right| ds
  \end{array}$$
By the Lipschitz continuity of $v^{\alpha}$ (uniform in $\alpha$)
that of $\psi$  and Proposition \ref{prop lip X}  we immediately
get:
\begin{equation}\label{lip-of-Z}
    \mathbb{E}\int_0^T |Z^{x,\alpha}_t-Z^{x',\alpha}_t|_{\Xi^*}^2 dt \leq c |x-x'|^2.
\end{equation}
for a suitable constant $c$ (that may depend on $T$).

Now we fix an arbitrary $T>0$ and, by a diagonal procedure  (using
separability of $E$) we  construct a subsequence
$(\alpha_n')\subset (\alpha_n)$ such that $\alpha_n' \searrow 0$
and
$$\mathbb{E}\int_0^T |Z^{x,\alpha_n'}_t-Z^{x',\alpha_m'}_t|_{\Xi^*}^2 dt \leq 2^{-n}
$$ for all $m\geq n$ and for all $x\in E$.
Consequently  $Z^{x,\alpha_n'}_t\rightarrow \overline{Z}^x_t$,
$\mathbb{P}$-a.s. for a.e. $t\in [0,T]$. Then we set:
$$\bar \zeta(x)=\left\{\begin{array}{ll} \lim_n \zeta^{\alpha_n'}(x),
& \hbox{ if the limit exists in }\Xi^*,\\
0,  & \hbox{ elsewhere.}\end{array}\right.$$ Since
$Z^{x,\alpha_n'}_t= \zeta^{\alpha_n'}(X^x_t)\rightarrow
\overline{Z}^{x}_t$   $\mathbb{P}$-a.s. for a.e. $t\in [0,T]$ we
immediately get that, for all $x\in E$, the process $X^x_t$
belongs   $\mathbb{P}$-a.s. for a.e. $t\in [0,T]$ to the set where
$\lim_n \zeta^{\alpha_n'}(x)$ exists  and consequently
 $\overline{Z}^{x}_t=\bar \zeta(X^x_t)$.
\end{proof}
\begin{remark}\begin{em} We notice that the solution we
have constructed above has the following ``linear growth''
property with respect to $X$: there exists $c>0$ such that,
$\mathbb{P}$-a.s.,
\begin{equation}\label{growt-of-Y}
|\overline{Y}^x_t|\leq c |X^x_t|  \hbox{ for all $t\geq 0$}.
\end{equation}
\end{em}
\end{remark}
If we require similar conditions then we immediately obtain
uniqueness of $\lambda$.
\begin{theorem}\label{th-uniq-lambda} Assume that,
in addition to Hypotheses \ref{general_hyp_forward}, \ref{hyp_W_A
F(W_A)} and \ref{hypothesisroyer}, Hypothesis
\ref{hyp-convol-determ} holds as well.  Moreover suppose that, for
some $x\in E$, the triple $(Y',Z',\lambda')$ verifies
$\mathbb{P}$-a.s. equation
 (\ref{EBSDE}) for all $0\leq t\leq T$,
where
 $Y'$ is a progressively measurable continuous process, $Z'$ is a process
 in $L_{\cal P, {\rm loc}}^2(\Omega;L^2(0,\infty;\Xi^*))$ and
 $\lambda'\in \mathbb{R}$.
 Finally assume that there exists $c_x>0$ (that may depend
 on $x$) such that
$\mathbb{P}$-a.s.
$$
 |Y'_t|\leq c_x (|X^x_t|+1) , \hbox{ for all $t\geq 0$}.
$$ Then $\lambda'=\bar \lambda$.
\end{theorem}
\begin{proof}
Let $\tilde \lambda=\lambda'-\lambda$, $\tilde
Y=Y'-\overline{Y}^x$, $\tilde Z=Z'-\overline{Z}^x$. By easy
computations:
$$\tilde \lambda=T^{-1}\left[\tilde Y_T-\tilde Y_0\right]+T^{-1}\int_0^T  \tilde Z_t \gamma_t dt
-T^{-1}\int_0^T  \tilde Z_t dW_t$$ where
$$\gamma_t:=\begin{cases}   \frac{\displaystyle \psi(X^{x}_t,Z'_t)-\psi(X^{x}_t,\overline{Z}^{x}_t)}{\displaystyle |Z'_t - \overline{Z}_t|_{\Xi^*}^2}\left(Z'_t - \overline{Z}_t \right)^*,& \hbox{ if } Z'_t \neq \overline{Z}_t, \\
0, & \hbox{  elsewhere   },
 \end{cases}
$$
is a bounded $\Xi$-valued progressively measurable process. By the
Girsanov Theorem  there exists a probability measure
$\mathbb{P}_{\gamma}$ under which $W^{\gamma}_t=-\int_0^t \gamma_s
ds+W_t$, $t\in [0,T]$, is a cylindrical Wiener process in $\Xi$.
Thus computing expectation with respect to $\mathbb{P}_{\gamma}$
we get
$$\tilde \lambda=T^{-1}\mathbb{E}^{\mathbb{P}_{\gamma}}
\left[\tilde Y_T-\tilde Y_0\right].$$ Consequently, taking into
account (\ref{growt-of-Y}),
 \begin{equation}\label{eq-proof-uniq-lambda}
|\tilde \lambda|\leq c T^{-1}\mathbb{E}^{\mathbb{P}_{\gamma}}
(|X^x_T|+1)+ c T^{-1}(|x|+1)
   \end{equation}
With respect to $W^{\gamma}$, $X^x$ is the mild solution of
$$
\left\{
\begin{array}{l}
dX^{x,\gamma}_t  =AX^{x,\gamma}_{t} dt+F( X^{x,\gamma}_{t} )
dt+GdW^{\gamma}_{t}+G\gamma_{t}\,dt ,
\quad t\geq 0 \\
X^{x,\gamma}_{0} =x\in E.
\end{array} \right.
$$
and by (\ref{rel-estimate-Xgamma}) we get
$\sup_{T>0}\mathbb{E}^{\mathbb{P}_{\gamma}}|X^x_T|<\infty$. So if
we let $T\rightarrow\infty$ in  (\ref{eq-proof-uniq-lambda}) we
conclude that $\tilde\lambda=0$.
\end{proof}

\begin{remark} \em
The solution to EBSDE (\ref{EBSDE}) is, in general,  not unique.
It is evident that the equation is invariant with respect to
addition of a constant to $Y$ but we can also construct an
arbitrary number of solutions that do not differ only by a
constant (even if we require them to be bounded). On the contrary
the solutions we construct are not Markovian.

Indeed, consider the equation:
\begin{equation}\label{eq:nouniqueness}
-dY_t=[\psi(Z_t)-\lambda]dt-Z_tdW_t.
\end{equation}
where $W$ is a standard brownian motion  and
$\psi:\mathbb{R}\rightarrow \mathbb{R}$ is differentiable bounded
and  has bounded derivative.

One solution is $Y=0;Z=0;\lambda=\psi(0)$ (without loss of
generality we can suppose that $\psi(0)=0$).

Let now  $\phi:\mathbb{R}\rightarrow \mathbb{R}$ be an arbitrary
differentiable function bounded and  with bounded derivative. The
following BSDE on $[t,T]$ admits a solution:
$$\left\{\begin{array}{rcl}
-dY_s^{x,t}&=&\psi(Z_s^{x,t})ds-Z_s^{x,t}dW_s,\\
Y_T^{x,t}&=&\phi(x+W_T-W_t).
\end{array}\right.$$
If we define $u(t,x)=Y_t^{x,t}$ then both $u$ and $\nabla u$ are
bounded.  Moreover if $\tilde{Y}_t=Y_t^{0,0}=u(t,W_t),\
\tilde{Z}_t=Z_t^{0,0}=\nabla u(t,W_t)$ then
$$\left\{\begin{array}{rcl}
-d\tilde{Y}_t&=&\psi(\tilde{Z}_t)dt-\tilde{Z}_tdW_t,\quad
t\in [0,T],\\
\tilde{Y}_T&=&\phi(W_T).
\end{array}\right.$$
 Then it is enough to extend with
$\tilde{Y}_t=\tilde{Y}_T,\ \tilde{Z}_t=0$ for $t>T$ to construct a
bounded solution to (\ref{eq:nouniqueness}).
\end{remark}
\begin{remark}\em The existence result in Theorem  \ref{main-EBSDE}
can be easily extended to the case of $\psi$ only satisfying
the conditions
$$ |\psi(x,z)
-\psi(x',z)|\le K_x|x-x'|,\quad  |\psi(x,0)|\le M,
\quad |\psi(x,z)|
\le K_z(1+|z|).
$$
Indeed we can construct a sequence $\{\psi_n : n\in \mathbb{N}\}$
of functions Lipschitz in $x$ and $z$ such that for all $x,x'\in
H$, $z \in \Xi^*$, $n\in \mathbb{N}$
$$ |\psi^n(x,z)
-\psi^n(x',z)|\le K'_x|x-x'|;\quad  |\psi^n(x,0)|\leq M';\quad
\lim_{n\rightarrow \infty}|\psi^n(x,z) -\psi(x,z)|=0.
$$
This can be done by projecting $x$ to the subspaces generated by
a basis in $\Xi^*$ and then regularizing by the standard
mollification techniques,  see \cite{FuTeBE}.
We know that if $(\bar Y^{x,n}, \bar Z^{x,n},\lambda_n)$ is the
solution of the EBSDE (\ref{EBSDE}) with $\psi$ replaced by
$\psi^n$ then $\bar Y^{x,n}_t=\bar v^n(X^x_t)$ with
$$ |\bar v^n(x)
-\bar v^n(x')|\le \dfrac{K'_x}{\eta}|x-x'|;\quad  \bar v^n(0)=0
;\quad |\lambda_n|\leq M'
$$
Thus we can assume (considering, if needed, a subsequence) that
$\bar v^n(x) \rightarrow \bar v(x)$ and $\lambda_n \rightarrow
\lambda$.
The rest of the proof is identical to the one of Theorem
\ref{main-EBSDE}.
\end{remark}

\section{Differentiability}

We are now interested in the differentiability of the
solution to the EBSDE (\ref{EBSDE}) with respect to $x$.

\begin{theorem}\label{th-diff} Assume that Hypotheses
\ref{general_hyp_forward} and
\ref{hypothesisroyer} hold. Moreover assume that $F$ is of class
${\cal G}^1(E,E)$ with $\nabla F$ bounded on bounded sets of $E$.
Finally assume that $\psi$ is of class ${\cal G}^1(E\times
\Xi^*,E)$. Then the function $\overline{v}$ defined in
(\ref{def-of-v}) is of class ${\cal G}^1(E,\mathbb{R})$.
\end{theorem}
\begin{proof} In \cite{masiero} it is proved that for arbitrary $T>0$ the map
$x\rightarrow X^x$ is of class $\mathcal{G}^1$ from $E$ to
$L^p_{\mathcal{P}}(\Omega,C([0,T],E))$. Moreover Proposition
\ref{prop lip X} ensures that for all $h\in E$,
\begin{equation}\label{proof-diff-estim-nabla-X}
 |\nabla X^x_t h|\leq e^{-\eta t}|h|,\quad \hbox{$\mathbb{P}$-a.s.,
 for all $t\in [0,T]$}.
\end{equation}
Under the previous conditions one can proceed exactly as
in Theorem 3.1 of \cite{HuTe} to
prove that for all $\alpha >0$ the map $v^{\alpha}$ is of class
$\mathcal{G}^1$.

$ $

Then we consider again
equation (\ref{bsderoyer}):
$$
{Y}^{x,\alpha}_t ={Y}^{x,\alpha}_T
+\int_t ^T(\psi(X^{x}_\sigma, Z^{x,\alpha}_\sigma)-\alpha
Y^{x,\alpha}_\sigma)d\sigma-\int_t ^T Z^{x,\alpha}_\sigma
dW_\sigma, \quad 0\le t \le T <\infty,
$$
we recall that ${Y}^{x,\alpha}_T={v}^{\alpha}(X^{x}_T)$,
 and apply again \cite{masiero} (see Proposition 4.2 there) and \cite{FuTe1}
 (see Proposition 5.2 there) to obtain that for all $\alpha >0 $
 the map $x\rightarrow Y^{x,\alpha}$ is of class $\mathcal{G}^1$
 from $E$ to $L^2_{\mathcal{P}}(\Omega,C([0,T],\mathbb{R}))$ and the map
$x\rightarrow Z^{x,\alpha}$  is of class $\mathcal{G}^1$ from $E$
to $L^2_{\mathcal{P}}(\Omega,L^2([0,T],\Xi^*))$. Moreover for all
$h\in E$ it holds (for all $t>0$ since $T$ was arbitrary)
$$
-d\nabla Y^{\alpha,x}_th=[\nabla_x\psi(X^x_t,Z_t^{\alpha,x})
\nabla X_t^xh+\nabla_z\psi(X^x_t,Z_t^{\alpha,x})\nabla
Z_t^{\alpha,x}h-\alpha\nabla Y^{\alpha,x}_th]dt
-\nabla Z^{\alpha,x}h dW_t.
$$
We also know that $|Y^{\alpha,x}_t|\le {M}/{\alpha}$. Now we set
$$U^{\alpha,x}_t=e^{\eta t}\nabla Y^{\alpha,x}_t h,
\quad V^{\alpha,x}=e^{\eta t}\nabla Z^{\alpha,x}_t h.$$ Then
$(U^{\alpha,x},V^{\alpha,x})$ satisfies the following BSDE:
\begin{eqnarray*}
-dU^{\alpha,x}_t&=&[e^{\eta t}\nabla_x\psi(X^x_t,Z_t^{\alpha,x})
\nabla X_t^x-(\alpha+\eta)U^{\alpha,x}_t +\nabla_z
\psi(X^x_t,Z_t^{\alpha,x}) V^{\alpha,x}_t]dt-V^{\alpha,x}_tdW_t.
\end{eqnarray*}
By (\ref{proof-diff-estim-nabla-X}) and the usual Girsanov
argument (recall the  $\nabla_x \psi$ and $\nabla_z \psi$ are
bounded),
$$|U^{\alpha,x}_t|\le \frac{c}{\alpha+\eta},\;
\forall t\geq 0,\; \hbox{$\mathbb P-$a.s. $\qquad$  i.e. } \qquad
|\nabla Y_t^{x,\alpha}|\le e^{-\eta t}\frac{c}{\alpha+\eta}.$$
Moreover, consider the limit equation, with unknown
$(U^{x},V^{x})$,
\begin{equation}\label{eq:limit}
-dU^x_t=[e^{\eta t}\nabla_x\psi(X^x_t,\bar Z_t^{x}) \nabla
X_t^x-\eta U^x_t+\nabla_z\psi(X^x_t,\bar Z_t^{x}) V^x]dt-V^xdW_t,
\end{equation}
which, since $|e^{\eta t}\nabla_x\psi \nabla_x X_t^x|$ is bounded,
has a unique  solution such that $U^x$ is bounded and $V^x$
belongs to $L_{\cal P, {\rm loc}}^2(\Omega;L^2(0,\infty;\Xi^*))$
(see \cite{bh} and \cite{royer}).

We know that for a suitable sequence $\alpha_n \searrow 0$,
$$\bar v^{\alpha}(x)= Y^{x,\alpha_n}_0-Y^{0,\alpha_n}_0\rightarrow \bar{Y}^x_0,$$
and we claim now that
$$  \nabla \bar v^{\alpha_n}(x)=\nabla Y_0^{x,\alpha_n}=U_0^{x,\alpha_n}
\rightarrow U_0^x.$$ To prove this we introduce the finite horizon
equations: for $t\in [0,N]$,
$$\begin{cases}
& -dU_t^{x,\alpha,N}=[e^{\eta t}\nabla_x\psi(X^x_t,Z_t^{x,\alpha})
 \nabla X_t^x-(\alpha+\eta)U_t^{x,\alpha,N}
 +\nabla_z \psi (X^x_t,Z_t^{x,\alpha}) V_t^{x,\alpha,N}]dt\\
& \qquad\qquad\qquad - V^{x,\alpha,N}_tdW_t,\\
& U_N^{x,\alpha,N}=0.
\end{cases}$$
$$\begin{cases}& -dU_t^{x,N}=[e^{\eta t}\nabla_x\psi(X^x_t,\bar Z_t^{x})
\nabla X_t^x-(\alpha+\eta)U_t^{x,N}
+\nabla_z \psi (X^x_t,\bar Z^{x}_t) V_t^{x,N}]dt-V^{x,N}_tdW_t,\\
& U_N^{x,N}=0.
\end{cases}$$
Since $\displaystyle \E\int_0^N |Z^{x,\alpha_n}_s-\bar Z^{x}_s|^2
ds\rightarrow 0$ it is easy to verify that, for all fixed $N>0$,
$U_0^{x,\alpha_n,N}\rightarrow U_0^{x,N}$.

On the other side a standard application of Girsanov Lemma gives
 see   \cite{HuTe},
$$|U_0^{x,\alpha_n,N}-U_0^{x,\alpha_n}|\le \frac{c}{\alpha_n+\eta}e^{-\eta N}, \qquad |U_0^{x,N}-U_0^{x}|\le \frac{c}{\eta}e^{-\eta N}.$$
for a suitable constant $c$.

Thus a standard argument implies $U_0^{x,\alpha_n}\rightarrow
U_0^{x}$. An identical argument also ensures continuity of
$U_0^{x}$ with respect to $x$ (also taking into account
\ref{lip-of-Z}). The proof is therefore completed.
\end{proof}

$ $

As  usual in the theory of markovian  BSDEs, the differentiability
property allows to identify the process  $\bar Z^x$ as a function
of the process $X^x$. To deal with our Banach space setting we
need to make the following extra assumption:

\begin{hypothesis}\label{Hyp-masiero}
There exists a Banach space $\Xi_0$, densely and continuously
embedded in $\Xi$, such that $G\, (\Xi_0) \subset \Xi$  and $G
:\Xi_0 \rightarrow E$ is continuous.
\end{hypothesis}

We note that this condition is satisfied in most applications. In
particular it is trivially true in the special case $E=H$ just by
taking $\Xi_0=\Xi$, since  $G$ is assumed to be a linear bounded
operator from $\Xi$ to $H$. The following is proved in
\cite[Theorem 3.17]{masiero}:

\begin{theorem}
 \label{theorem-identif-Z}
Assume that Hypotheses \ref{general_hyp_forward},
 \ref{hypothesisroyer} and \ref{Hyp-masiero} hold.
Moreover assume that $F$ is of class ${\cal G}^1(E,E)$ with
$\nabla F$ bounded on bounded subsets of $E$ and $\psi$ is of
class ${\cal G}^1(E\times \Xi^*,E)$. Then $\bar Z^x_t=\nabla \bar
v(X^x_t)G$, $\mathbb{P}$-a.s. for a.e. $t\geq 0$.
\end{theorem}
\begin{remark} \label{precision}\begin{em}
We notice that $\nabla \bar v(x)G\xi$ is only defined for $\xi\in
\Xi_0$ in general, and the conclusion of
Theorem \ref{theorem-identif-Z} should be stated more precisely
as follows: for $\xi\in
\Xi_0$ the equality $Z^x_t\xi=\nabla \bar v(X^x_t)G\xi$
holds $\mathbb{P}$-a.s. for almost every $t\geq 0$. However,
since $\bar Z^x$
is a process with values in $\Xi^*$, and more specifically
a process in  $
L^2_{\mathcal{P}}(\Omega,L^2([0,T],\Xi^*))$, it follows that
$\P$-a.s. and
for almost every
$t$  the
operator $\xi \rightarrow \nabla \bar v(X^x_t)G\xi$ can be
extended to a bounded linear operator defined on the whole $\Xi$.
Equivalently,
for almost every
$t$ and for almost all $x\in E$  (with respect to the law of $X_t$)
the linear
operator $\xi \rightarrow \nabla \bar v(x)G\xi$ can be
extended to a bounded linear operator defined on the whole $\Xi$
(see also Remark 3.18 in \cite{masiero}).
\end{em}
\end{remark}
\begin{remark}\label{boundedpsibar}
\begin{em} The above representation together with the fact that
$\bar v$ is Lipschitz with Lipschitz constant $K_x\eta^{-1}$
immediately implies that, if $F$ is of class ${\cal G}^1(E,E)$ and
$\psi$ is of class ${\cal G}^1(E\times \Xi^*,E)$, then $|\bar
{Z}^x_t|_{\Xi_0^*}\leq K_x\eta^{-1} |G|_{L(\Xi_0,E)}$ for all $x\in
E$, $\mathbb{P}$-a.s. for almost every $t\geq 0$. Consequently we
can construct $\bar \zeta$ in Theorem \ref{main-EBSDE}
 in such a way that it is bounded in the
$\Xi_0^*$ norm by $K_x\eta^{-1} |G|_{L(\Xi_0,E)}$.

Once this is proved we can extend the result to the case in which
$\psi$ is no longer differentiable but only Lipschitz, namely
we can prove than even in this case the process $\bar
{Z}^x$ is bounded. Indeed if we
consider a sequence $\{\psi_n : n\in \mathbb{N}\}$ of functions of
class ${\cal G}^1(E\times \Xi^*,E)$ such that for all $x,x'\in H$,
$z,z'\in \Xi^*$, $n\in \mathbb{N}$,
$$ |\psi_n(x,z)
-\psi_n(x',z')|\le K_x|x-x'|+ K_z |z-z'|;\quad \lim_{n\rightarrow
\infty}|\psi_n(x,z) -\psi(x,z)|=0.
$$
We know that if $(\bar Y^{x,n}, \bar Z^{x,n},\lambda_n)$ is the
solution of the EBSDE (\ref{EBSDE}) with $\psi$ replaced by
$\psi_n$ then $|\bar {Z}^{x,n}_t|_{\Xi_0^*}\leq K_x\eta^{-1}
|G|_{L(\Xi_0,E)}$. Then as we did above we can show (showing that the
corresponding equations with monotonic generator converge
uniformly in $\alpha$) that $\mathbb{E}\int_0^T|\bar {Z}^{x,n}_t
-\bar {Z}^{x}_t|_{\Xi_0^*}^2dt\rightarrow 0$ and the claim
follows.

We also notice that by the same argument we also have $ |\bar
\zeta^{\alpha}(x)|_{\Xi_0^*}\leq K_x\eta^{-1} |G|_{L(\Xi_0,E)}$,
$\forall \alpha>0$.
\end{em}
\end{remark}
Now we introduce  the Kolmogorov semigroup corresponding to $X$:
 for  measurable and bounded  $\phi:
E\rightarrow \mathbb{R}$ we define
\begin{equation}\label{def-of-p}
P_t[\phi](x)=\mathbb{E}\, \phi(X^x_t)\qquad t\ge 0,\, x\in E.
\end{equation}
\begin{definition}\label{strongly-feller}
The semigroup $(P_t)_{t\geq 0}$ is called strongly Feller if for
all $t>0$ there exists $k_t$ such that for all measurable and
bounded $\phi: E\rightarrow \mathbb{R}$,
$$|
P_t[\phi](x)- P_t[\phi](x')|\leq k_t \Vert\phi\Vert_0 |x-x'|,
\qquad x,x'\in E,
$$
where $\Vert\phi\Vert_0=\sup_{x\in E}|\phi(x)|$.
\end{definition}

This terminology is somewhat different from the classical one
(namely, that $P_t$ maps  measurable  bounded functions into
continuous ones,
 for all
$t>0$), but it will be convenient for us.

\begin{definition}\label{gen-diss} We say that $F$
is genuinely dissipative if there exist $\epsilon>0$ and $c>0$
such that, for all $x,x'\in E$, there exists $z^*\in \partial
|x-x'|$ such that $<z^*,F(x)-F(x')>_{E^*,E}\leq c
|x-x'|^{1+\epsilon}$.
\end{definition}

\begin{lemma}\label{lemma-SF-dissip}
Assume that Hypotheses \ref{general_hyp_forward} and \ref{hyp_W_A
F(W_A)} hold.
If the Kolmogorov
semigroup $(P_t)$ is strongly Feller then for all bounded
measurable $\phi: E\rightarrow\mathbb{R}$,
$$\left|P_t[\phi](x)-\int_E \phi(x)\mu (dx)\right|
\leq c e^{-\eta (t/4)}(1+|x|)\Vert\phi\Vert_0.$$
If in addition $F$ is genuinely dissipative then
$$\left|P_t[\phi](x)-\int_E \phi(x)\mu (dx)\right|
\leq c e^{-\eta (t/4)}\Vert\phi\Vert_0.$$
\end{lemma}
\begin{proof} We fix $\epsilon >0$. For $t>2$ we have,
by Theorem \ref{ergodicity},
$$
\begin{array}{r}\displaystyle\left|P_t[\phi](x)-\int_E \phi(x)\mu (dx)\right|=
\left|P_{t-1}[P_1[\phi]](x)-\int_E P_{1}[\phi](x)\mu (dx)\right|
\leq C(1+|x|)
e^{-\eta t /4} \Vert P_{1}[\phi]\Vert_{\hbox{lip}}\\
\displaystyle \leq C(1+|x|) e^{-\eta t /4} k_{1}\Vert\phi\Vert_0,
\end{array}$$
and the first claim follows since $\left|P_t[\phi](x)-\int_E
\phi(x)\mu (dx)\right|\leq 2  \Vert\phi\Vert_0$.

If now $F$ is genuinely dissipative then in \cite{DP2}, Theorem
6.4.1 it is shown that
$$\left|\mathbb{E}\phi(X^x_t)-\int_E \phi\, d\mu \right|\leq
C e^{-\eta t /2} \Vert\phi\Vert_{\hbox{lip}}$$ and the second
claim follows by the same argument.
\end{proof}

We are now able to state and prove two corollaries
of Theorems  \ref{th-diff}  and \ref{theorem-identif-Z}.

\begin{corollary}\label{characterization of lambda}
Assume that Hypotheses \ref{general_hyp_forward}, \ref{hyp_W_A
F(W_A)}, \ref{hypothesisroyer} and \ref{Hyp-masiero} hold.
Moreover assume that $F$ is of class $\mathcal{G}^1$ with $\nabla
F$ bounded on bounded subsets of $E$, and that $\psi$ is bounded
on each set $E\times B$, where $B$ is any ball of  $\Xi_0^*$.
Finally assume that the Kolmogorov semigroup $(P_t)$ is strongly
Feller.

Then the following holds:
$$\lambda=\int_E \psi(x,\bar \zeta(x))\mu (dx),$$
where $\mu$ is the unique invariant measure of $X$.
\end{corollary}
\begin{proof} First notice that $\overline{\psi}:=
\psi(\,\cdot\, , \bar\zeta(\,\cdot\,))$ is bounded, by
Remark \ref{boundedpsibar}.
 Then
$$T^{-1}\mathbb{E}[\bar Y ^x_0-\bar Y ^x_T]=
T^{-1}\E \int_0^T\left  (\psi(X^x_t,\bar \zeta( X^x_t))- \int_E
\bar \phi\, d\mu \right)dt+ \left(\int_E \bar \phi\, d\mu
-\lambda\right).$$ We know that $T^{-1}\mathbb{E}[\bar Y ^x_0-\bar
Y ^x_T]\rightarrow 0$, by the argument
in Theorem \ref{th-uniq-lambda}.
Moreover by the first conclusion of Lemma
\ref{lemma-SF-dissip}
$$ T^{-1}\E \int_0^T\left  (\psi(X^x_t,\bar \zeta( X^x_t))-
\int_E \bar \phi\, d\mu \right)dt \rightarrow 0,$$ and the claim
follows. \end{proof}

\begin{corollary}\label{boundedness of v}
In addition to the assumptions of Corollary \ref{characterization
of lambda} suppose that $F$ is genuinely dissipative. Then $\bar
v$ is bounded.
\end{corollary}
\begin{proof}
Let $(Y^{x,\alpha},Z^{x,\alpha})$ be the solution of
(\ref{bsderoyer}). We know that $Y^{x,\alpha}_t=v^{\alpha}(X^x_t)$
and $Z^{x,\alpha}_t= \zeta^{\alpha}(X^x_t)$  with $v^{\alpha}$
Lipschitz uniformly with respect to $\alpha$ and $\zeta^{\alpha}$
bounded in $\Xi^*$ uniformly with respect to $\alpha$. Let
$\psi^{\alpha}=\psi(\,\cdot\,,\bar \zeta^{\alpha}(\,\cdot\,))$.
Under the present assumptions we conclude that also the maps
$\psi^{\alpha}$ as well are bounded in $\Xi^*$ uniformly with
respect to $\alpha$.

Computing $d (e^{-\alpha t} \bar Y^{x\alpha}_t)$ we obtain,
$$Y^{x,\alpha}_0=\mathbb{E} e^{-\alpha T} Y^{x,\alpha}_T+
\mathbb{E} \int_0^T  e^{-\alpha t} \psi^{\alpha} (X^x_t)dt,$$ and
for $T\rightarrow\infty$,
$$Y^{x,\alpha}_0=
\mathbb{E} \int_0^\infty  e^{-\alpha t} \psi^{\alpha} (X^x_t)dt.$$
Subtracting  to both sides $\alpha^{-1}\int_E
\psi^{\alpha}(x)\mu(dx)$ we obtain
$$\left|Y^{x,\alpha}_0-\alpha^{-1}\int_E  \psi^{\alpha}(x)\mu(dx)\right|=
\left| \int_0^\infty  e^{-\alpha t} \left[P_t[\psi^{\alpha}]
(x)-\int_E  \psi^{\alpha}(x)\mu(dx)\right]dt\right|\leq 4c
\eta^{-1} \Vert \psi^\alpha\Vert_0 $$ where the last inequality
comes from the second conclusion of Lemma \ref{lemma-SF-dissip}.

Thus  $\left|Y^{x,\alpha}_0-Y^{0,\alpha}_0\right| \leq 8 c
\eta^{-1} \Vert \psi^\alpha\Vert_0 $ and the claim follows since
by construction $Y^{x,\alpha}_0-Y^{0,\alpha}_0  \rightarrow \bar v
(x)$.
\end{proof}
\section{Ergodic Hamilton-Jacobi-Bellman equations}
We briefly show here that if $\bar Y_0^x=\bar v(x)$ is of class
${\cal G}^1$ then the couple $(v,\lambda)$  is a mild solution of
the following ``ergodic''  Hamilton-Jacobi-Bellman equation:
\begin{equation}
\mathcal{L}v(x)
+\psi\left( x,\nabla v(x)  G\right)  = \lambda, \quad x\in E,  \label{hjb}%
\end{equation}
Where  linear operator $\mathcal{L}$ is formally defined by
\[
\mathcal{L}f\left(  x\right)  =\frac{1}{2}Trace\left(
GG^{\ast}\nabla ^{2}f\left(  x\right)  \right)  +\langle Ax,\nabla
f\left(  x\right) \rangle_{E,E^{\ast}}+\langle F\left(  x\right)
,\nabla f\left(  x\right) \rangle_{E,E^{\ast}},
\]
We notice that we can define the transition semigroup
 $(P_t)_{t\geq 0}$ corresponding to $X$ by the formula (\ref{def-of-p})
for all measurable functions $\phi:E\to\mathbb{ R}$ having
polynomial growth, and we notice that $\mathcal{L}$ is the formal
generator of $(P_t)_{t\geq 0}$.

 Since we are dealing with an elliptic equation it is natural to consider
$(v,\lambda)$ as a mild solution of equation (\ref{hjb}) if and
only if, for arbitrary $T>0$, $v(x)$ coincides with the  mild
 solution $u(t,x)$ of the corresponding parabolic equation
 having $v$ as a terminal condition:
\begin{equation}\left\{
\begin{array}{l}
 \dfrac{\partial u(t,x)}{\partial t}+\mathcal{L}u\left(  t,x\right)
+\psi\left(  x,\nabla u\left(  t,x\right)  G\right)
 -\lambda=0, \quad t\in [0,T],\; x\in E,  \\ \\
u(T,x)=v(x), \quad  x\in E.
 \end{array}\right. \label{hjb-parab}
\end{equation}
Thus we are led to the following definition (see also
\cite{FuTe-ell}):
\begin{definition}
\label{defsolmildkolmo} A pair $(v,\lambda)$ ($v: E\rightarrow
\mathbb{R}$ and $\lambda\in \mathbb{R}$) is a mild solution of the
Hamilton-Jacobi-Bellman equation (\ref{hjb}) if the following are
satisfied:

\begin{enumerate}
\item $v\in\mathcal{G}^{1}\left(  E,\R \right)  $;

\item  there exists $C>0$ such that $\left|  \nabla v\left(
x\right)
h\right|  \leq C\left|  h\right|  _{E}\left(  1+\left|  x\right|  _{E}%
^{k}\right)  $ for every  $x,h\in E$ and some positive integer
$k$;

\item for $0\le t\le T$  and $x\in E$,
\begin{equation}
v(x)=P_{T-t}\left[  v\right]  \left(  x\right)
+\int_{t}^{T}\left(P_{s -t }\left[  \psi(\cdot,\nabla v\left(
\cdot\right)  G)\right]  \left( x\right) -\lambda \right)   \,ds.
\label{mild sol hjb}
\end{equation}

\end{enumerate}
\end{definition}

In the right-hand side of (\ref{mild sol hjb}) we notice
occurrence of the term $\nabla v\left( \cdot\right)  G$, which is
not well defined as a function $E\to\Xi^*$, since $G$ is not
required to map $\Xi$ into $E$.
The situation is similar to Remark \ref{precision}.
In general,
 for $x \in E$, $\nabla \bar
v(x)G\xi$ is only defined for $\xi\in \Xi_0$.
In (\ref{mild sol hjb}) it is implicitly required
that, $\P$-a.s. and
for almost every
$t$,  the
operator $\xi \rightarrow \nabla \bar v(X^x_t)G\xi$ can be
extended to a bounded linear operator defined on the whole $\Xi$.
Noting that
$$
P_{t }\left[  \psi(\cdot,\nabla v\left( \cdot\right)  G)\right]
\left( x\right) = \E \, \psi(X^x_{t},\nabla v\left( X^x_{t}\right)
G)
$$
the equation (\ref{mild sol hjb}) is now meaningful.

Using the results for the parabolic case,  see \cite{masiero}, we
get existence of the mild solution of equation (\ref{hjb})
whenever we have proved that the function
$\bar v$ in Theorem \ref{main-EBSDE} is differentiable.

\begin{theorem}\label{th-EHJB}
Assume that Hypotheses \ref{general_hyp_forward},
\ref{hypothesisroyer} and \ref{Hyp-masiero} hold.
Moreover assume that $F$ is of class ${\cal G}^1(E,E)$ with
$\nabla F$ bounded on bounded subsets of $E$ and $\psi$ is of
class ${\cal G}^1(E\times \Xi^*,E)$.

Then $(\bar v,  \bar\lambda)$ is a   mild solution of the
Hamilton-Jacobi-Bellman equation (\ref{hjb}).

Conversely, if $(v,\lambda)$ is a   mild solution of
 (\ref{hjb}) then, setting  $ Y^x_t=
v(X^x_t)$ and ${Z}^{x}_t=
\nabla v( X^x_t)  G$,
the triple
 $( Y^x, Z^x,  \lambda)$ is a solution of
  the EBSDE
 (\ref{EBSDE}).

\end{theorem}

\section{Optimal ergodic control}
\label{optcontr}

Assume that Hypothesis \ref{general_hyp_forward} holds and let
$X^x$ denote the solution to equation (\ref{sde}).
 Let $U$ be a separable
 metric space. We define a control $u$ as an
$(\calf_t)$-progressively measurable $U$-valued process.  The cost
 corresponding to a given control
is defined in the following way. We assume that the functions
$R:U\rightarrow \Xi^*$ and $L:E\times U \rightarrow \R$ are
measurable and satisfy, for some constant $c>0$,
\begin{equation}\label{condcosto}
|R(u)|\leq c,\quad |L(x,u)|\leq c, \quad |L(x,u)-L(x',u)|\leq
c\,|x-x'|,\qquad u\in U,\,x,x'\in E.
\end{equation}
Given an arbitrary control $u$ and $T>0$, we introduce the
Girsanov density
$$ \rho_T^u=\exp\left(\int_0^T R(u_s)dW_s
-\frac{1}{2}\int_0^T |R(u_s)|_{\Xi^*}^2 ds\right)$$ and the
probability $\mathbb P_T^u=\rho_T^u\mathbb P$ on $\calf_T$. The
ergodic cost  corresponding to $u$ and the starting point $x\in E$
is
\begin{equation}\label{def-ergodic-cost}
  J(x,u)=\limsup_{T\rightarrow\infty}\frac{1}{T} \mathbb
E^{u,T}\int_0^T L(X_s^x,u_s)ds,
\end{equation}
where $\mathbb E^{u,T}$  denotes expectation with respect to
$\mathbb P_T^u$. We notice that $W_t^u=W_t-\int_0^t R(u_s)ds$ is a
Wiener process on $[0,T]$ under $\mathbb P^u$ and that
$$dX_t^x=(AX_t^x+F(X_t^x))dt+G(dW_t^u+R(u_t)dt),
\quad t\in [0,T]$$ and this justifies our formulation of the
control problem. Our purpose is to minimize the cost over all
controls.

 To this purpose we first define the Hamiltonian in the
usual way
\begin{equation}\label{defhamiton}
\psi(x,z)=\inf_{u\in U}\{L(x,u)+z R(u)\},\qquad x\in E,\,z\in
\Xi^*,
\end{equation}
and we remark that  if, for all $ x,z$, the  infimum is attained
in (\ref{defhamiton}) then there exists a measurable function
$\gamma:E\times \Xi^*\rightarrow U$ such that
$$\psi(x,z)=l(x,\gamma(x,z))+z R(\gamma(x,z)).$$
This follows from an application of Theorem 4 of \cite{McS-War}.

We notice that under  the present assumptions $\psi$ is a
Lipschitz function and $\psi(\cdot,0)$ is bounded (here the fact
that $R$ depends only on $u$ is used). So if we assume Hypotheses
\ref{general_hyp_forward} and \ref{hyp_W_A F(W_A)} then in Theorem
\ref{main-EBSDE} we have constructed, for every $x\in E$, a triple
\begin{equation}\label{richiamoebsde}
(\bar Y^x,\bar Z^x, \bar \lambda)= (\bar v (X^x),\bar \zeta(X^x),
\bar \lambda)
\end{equation} solution to
 the EBSDE
 (\ref{EBSDE}).

\begin{theorem}\label{Th-main-control}
Assume that Hypotheses \ref{general_hyp_forward}, \ref{hyp_W_A
F(W_A)} and  \ref{hyp-convol-determ} hold, and that
(\ref{condcosto}) holds as well.

Moreover suppose that, for some $x\in E$, a triple $(Y,Z,\lambda)$
verifies $\mathbb{P}$-a.s. equation
 (\ref{EBSDE}) for all $0\leq t\leq T$,
where
 $Y$ is a progressively measurable continuous process, $Z$ is a process
 in $L_{\cal P, {\rm loc}}^2(\Omega;L^2(0,\infty;\Xi^*))$ and
 $\lambda\in \mathbb{R}$.
 Finally assume that there exists $c_x>0$ (that may depend
 on $x$) such that
$\mathbb{P}$-a.s.
$$
 |Y_t|\leq c_x (|X^x_t|+1) , \hbox{ for all $t\geq 0$}.
$$

Then the following holds:
\begin{enumerate}
 \item[(i)] For arbitrary control
 $u$ we have $J(x,u)\ge \lambda=\bar\lambda,$
and the equality holds if and only if $L(X_t^x,u_t)+Z_t
R(u_t)=\psi(X_t^x,Z_t)$, $\P$-a.s. for almost every $t$.

\item[(ii)] If the  infimum is attained in (\ref{defhamiton}) then
the control $\bar u_t=\gamma(X_t^x,Z_t)$ verifies $J(x,\bar u)=
\bar\lambda.$
\end{enumerate}

In particular, for the solution (\ref{richiamoebsde}) mentioned
above, we have:
\begin{enumerate}
 \item[(iii)] For arbitrary control
 $u$ we have $J(x,u)=\bar\lambda$ if and only if
$L(X_t^x,u_t)+\bar\zeta (X_t^x) R(u_t)=\psi(X_t^x,\bar \zeta
(X_t^x))$, $\P$-a.s. for almost every $t$. \item[(iv)] If the
infimum is attained in (\ref{defhamiton}) then the control $\bar
u_t=\gamma(X_t^x,\bar\zeta (X_t^x))$ verifies $J(x,\bar u)=
\bar\lambda.$
\end{enumerate}

\end{theorem}

\begin{remark}\em
\begin{enumerate}
 \item
The equality $\lambda=\bar\lambda$ clearly follows from Theorem
\ref{th-uniq-lambda}. \item Points $(iii)$ and $(iv)$ are
immediate consequences of $(i)$ and $(ii)$. \item The conclusion
of point $(iv)$ is that there exists an optimal control in
feedback form, with the optimal feedback given by the function
$x\mapsto \gamma(x,\bar\zeta (x))$. \item Under the conditions of
Theorem \ref{th-EHJB}, the pair $(\bar v, \bar \lambda)$ occurring
in (\ref{richiamoebsde}) is a mild solution of the
Hamilton-Jacobi-Bellman equation (\ref{hjb}). \item It follows
from the proof below that if $\limsup$ is changed into $\liminf$
in the definition (\ref{def-ergodic-cost}) of the cost, then the
same conclusions hold, with the obvious modifications, and the
optimal value is given by $\bar\lambda$ in both cases.
\end{enumerate}
\end{remark}

\begin{proof}
 As $(Y,{Z}, \bar\lambda)$ is a solution of the
ergodic BSDE, we have
\begin{eqnarray*}
-d{Y}_t&=&[\psi(X_t^x,{Z}_t)-\bar\lambda]dt-{Z}_tdW_t\\
&=&[\psi(X_t^x,{Z}_t)- \bar\lambda]dt-{Z}_tdW_t^u-{Z}_t R(u_t)dt,
\end{eqnarray*}
from which we deduce that
\begin{eqnarray*}
\bar\lambda&=&\frac{1}{T}\mathbb E^{u,T}[Y_T-Y_0]
+\mathbb E^{u,T}\frac{1}{T}\int_0^T[\psi(X_t^x,{Z}_t)-{Z}_t r(u_t)-L(X_t^x,{Z}_t)]dt\\
& &+\frac{1}{T}\mathbb E^{u,T}\int_0^T L(X_t^x,{Z}_t)dt.
\end{eqnarray*}

Thus
$$\frac{1}{T}\mathbb E^{u,T}\int_0^T L(X_t^x,{Z}_t)dt\ge
 \frac{1}{T}\mathbb E^{u,T}[Y_0-Y_T]+\bar\lambda.$$
But by (\ref{rel-estimate-Xgamma}) we have
$$|\mathbb E^{u,T} Y_T|\le c\mathbb E^{u,T}(|X_T^x|+1)\le c(1+|x|).$$
Consequently $T^{-1}\mathbb E^{u,T}[Y_0-Y_T]\rightarrow 0,$ and
$$\limsup_{T\rightarrow\infty } \frac{1}{T}\mathbb E^{u,T}\int_0^T L(X_t^x,{Z}_t)dt
\ge \bar\lambda.$$

Similarly, if $L(X_t^x,u_t)+ Z_t R(u_t)=\psi(X_t^x,Z_t)$,
$$\frac{1}{T}\mathbb E^{u,T}\int_0^T L(X_t^x,{Z}_t)dt=
\frac{1}{T}\mathbb E^{u,T}[Y_0-Y_T]+\bar\lambda,$$ and the claim
holds.
\end{proof}

\section{Uniqueness}\label{sec-uniq}
We wish now to adapt the argument in \cite{GoMa} in order to
obtain uniqueness of markovian solutions to the EBSDE. This will
be done by a control thoretic interpretation the requires that the
Markov process related to the state equation with continuous
feedback enjoys recurrence properties. In this section we assume
\begin{equation}\label{addizionali}
E=H \qquad\hbox{ and }\qquad F  \hbox{ is bounded.}
\end{equation}

\noindent We recall here a result due to \cite{seid} on recurrence
of solution to SDEs.
\begin{theorem}\label{th-rec-seidler}
Consider
\begin{equation}\label{eq:u}
d{X}_t=(A{X}_t+g({X}_t))dt+GdW_t.
\end{equation}
where $g: H \rightarrow H$ is bounded and weakly continuous (that
if $x\rightarrow\<\xi,g(x)\>$ is continuous for all $\xi\in H$).
Let
$$Q_t=\int_0^t e^{sA}GG^*e^{sA^*}ds.$$
and assume the following
\begin{enumerate}
 \item $\sup_{t\ge 0} \hbox{Trace}\,(Q_t)<\infty$;
\item $Q_t$ is injective for $t>0$; \item  $ e^{t A}(H)\subset
(Q_t)^{1/2}(H)$ for $t>0$; \item $\int_0^t
|Q_s^{-1/2}e^{sA}|ds<\infty$ for $t>0$; \item there exists
$\beta>0$ such that $\int_0^t s^{-\beta}\,
\hbox{Trace}\,(S(s)S(s)^*)\, ds<\infty$
  for $t>0$.
\end{enumerate}
Then,  for all $T>0$, equation (\ref{eq:u}) admits a martingale
solution on $[0,T]$, unique in law. The associated transition
probabilities $P(t,x,T,\cdot)$ on $H$ ($0\le t\le T, x\in H$)
identify a   recurrent  Markov process on $[0,\infty)$.
\end{theorem}

Consider now the ergodic control problem with state equation:
$$d{X}^{x,u}_t=(A{X}^{x,u}_t+F({X}^{x,u}_t)+GR(u_t))dt+GdW_t, \ X_0^{x,u}=x,$$
and cost
$$\limsup_{T\to\infty}\frac{1}{T}\,
\mathbb E\int_0^T l(X_s,u_s)ds$$ where  $R:U\rightarrow
\Xi$ is continuous  and bounded.

We restrict ourselves to the class of controls given by continuous
feedbacks, i.e. given arbitrary
 continuous $u: H\rightarrow U$ (called feedback) we define the
 corresponding trajectory as the solution of
$$d{X}^{x,u}_t=(A{X}^{x,u}_t+F({X}^{x,u}_t))dt+G(R(u(X_t^{x,u}))dt+dW_t),
\ X_0^{u,x}=x.$$
We notice that for all $T>0$ there exists a weak solution $X^{x,u}$ of
this equation, and it is unique in law.

$ $

 We set as usual
$$\psi(x,z)=\inf_{u\in U}\{L(x,u)+zR(u)\},$$
and assume that $\psi$ is continuous and there exists a continuous
$\gamma:H\times\Xi\rightarrow U$ such that
$$\psi(x,z)=L(x,\gamma(x,z))+zR(\gamma(x,z)).$$

\begin{theorem}\label{th-uniqueness}
Suppose (\ref{addizionali})
and suppose that the assumptions of Theorem
\ref{th-rec-seidler} hold.
Let $(v,\zeta,\lambda)$ with $v:H\rightarrow \mathbb{R}$  continuous,
$\zeta:H\rightarrow \mathbb{R}$  continuous, and $\lambda$ a real number
satisfy the following conditions:

\begin{enumerate}
 \item $|v(x)|\le c|x|$;
\item
for an arbitrary  filtered probability space with
a Wiener process
$(\hat{\Omega},\hat{\mathcal{F}},
\{\hat{\mathcal{F}}_t\}_{t>0},\hat{\mathbb{P}},\{
\hat{W}_t\}_{t>0})$ and
for any solution of
$$d\hat{X}_t=(A\hat{X}_t+F(\hat{X}_t))dt+Gd\hat{W}_t,\qquad t\in [0,T],$$
setting  $Y_t=v(\hat{X}_t),\
Z_t=\zeta(\hat{X}_t)$, we have
$$-dY_t=[\psi(\hat{X}_t,Z_t)-\lambda]dt-Z_tdW_t\quad t\in [0,T].$$
\end{enumerate}
Let
$$\tau_r^T=\inf \{s\in [0,T]:|X_s^{u,x}|< r\},$$
with the convention $\tau_r^T=T$ if the indicated set
is empty,
and
$$J(x,u)=\limsup_{r\rightarrow 0}\limsup_{T\rightarrow \infty}
\mathbb E\int_0^{\tau_r^T} [\psi(X_s^{x,u},u(X_s^{x,u}))-\lambda]ds.$$
Then
$$v(x)=\inf_u J(x,u),$$
where the infimum (that is a minimum) is taken over all continuous
feedbacks $u$.
\end{theorem}
\begin{proof} Let $u:H\to U$ be continuous.
We notice that $X^{x,u}$ solves on $[0,T]$:
$$dX_t^{x,u}=(AX_t^{x,u}+F(X_t^{x,u}))ds+Gd\tilde{W}_t^u,\ t\in [0,T],$$
where $\tilde{W}_t=\int_0^t R(u(X_r^{x,u})dr+W_t$ is a Wiener
process on $[0,T]$ under a suitable probability $\hat{\mathbb{P}}^{u,T}$.

Therefore $Y_t=v(X_t^{x,u})$, $Z_t=\zeta(X_t^{x,u})$ satisfy:
$$
-dY_t=[\psi(X_t^{x,u},u(X_t^{x,u}))-\lambda]dt-Z_t
R(u(X_t^{x,u}))]dt-Z_tdW_t.$$ Integrating in $[0,\tau_r^T]$ we get
$$v(x)=\mathbb E(v(X_{\tau_r^T}^{x,u}))+\mathbb E\int_0^{\tau_r^T}
[\psi(X_s^{u,x},u(X_s^{x,u}))-\lambda-Z_s R(X_s^{x,u})]ds.$$ Thus,
\begin{equation}\label{eq:x}
v(x)\le\mathbb E(v(X_{\tau_r^T}^{x,u}))+\mathbb E\int_0^{\tau_r^T}
[L(X_s^{u,x},u(X_s^{x,u}))-\lambda]ds.
\end{equation}
Now
\begin{eqnarray*}
|\mathbb E(v(X_{\tau_r^T}^{x,u}))|\le c\mathbb
E|X^{x,u}_{\tau_r^T}|&\le& c r+ (\mathbb
E(|X_T^{x,u}|^2))^{1/2}(\mathbb P(\tau_r^T=r))^{1/2}\\ &\le & c r+
c (\mathbb P(\tau_r^T=r))^{1/2}
\end{eqnarray*}
Notice that $\mathbb P(\tau_r^T=r)=\tilde{\mathbb P}(\inf_{t\in
[0,T]}|\tilde{X}_t|\geq r),$ where $\tilde{X}$ is the Markov
process on the whole $[0,+\infty)$ corresponding to the equation
(\ref{eq:u}) with $g=F(\cdot)+GR(u(\cdot))$.

$ $

\noindent Since $\tilde{X}$ is recurrent, for all $ r>0$ it holds
$\tilde{\mathbb P}(\inf_{t\in [0,T]}|\tilde{X}_t|>r)\rightarrow 0$
as $T\rightarrow \infty.$ Thus
$$\limsup_{r\rightarrow 0}\limsup_{T\rightarrow \infty}|
\mathbb E(v(X_{\tau_r^T}^{x,u}))|\rightarrow 0.$$
Hence,
$$v(x)\le \limsup_{r\rightarrow 0}\limsup_{T\rightarrow \infty}
\mathbb E\int_0^{\tau_r^T}
[l(X_s^{x,u},u(X_s^{x,u}))-\lambda]ds.$$ The proof is completed
noticing that if  $u$ is chosen as ${u}(x)=\gamma(x,\zeta(x))$
then the  above
inequality becomes an equality.
\end{proof}

This result combines with Theorems \ref{th-uniq-lambda}
and \ref{th-EHJB}
to give the following

\begin{corollary}\label{HJB-uniqueness}
Suppose  that all the assumptions of
Theorems \ref{th-uniq-lambda}, \ref{th-EHJB} and \ref{th-uniqueness} hold.
Then  $(\bar v,  \bar\lambda)$ is the unique   mild solution of the
Hamilton-Jacobi-Bellman equation (\ref{hjb}) satisfying
$|\bar v (x)|\le c|x|$.
\end{corollary}

\section{Application to ergodic control of a semilinear heat equation}
\label{section-heat-eq}

In this section we show how our results can be applied to perform
the synthesis of the ergodic optimal control when the state
equation is a semilinear heat equation with additive noise.  More
precisely, we treat a stochastic heat equation in space dimension
one, with a dissipative nonlinear term and with control and noise
acting on a subinterval. We consider homogeneous Dirichlet
boundary conditions.

\noindent In $\left(  \Omega,\mathcal{F},\mathbb{P}\right)  $ with
a filtration $\left(  \mathcal{F}_{t}\right) _{t\geq0}$ satisfying
the usual conditions, we consider, for $t \in\left[ 0,T\right] $
and $\xi\in\left[  0,1\right]  $, the following equation
\begin{equation}
\left\{
\begin{array}
[c]{l} d_{t }X^{u}\left(  t ,\xi\right)  =\left[
\frac{\partial^{2}}{\partial \xi^{2}}X^{u}\left(  t ,\xi\right)
+f\left(  \xi,X^{u}\left(  t  ,\xi\right)  \right)
+\chi_{[a,b]}(\xi) u\left( t ,\xi\right)  \right] dt
+\chi_{[a,b]}(\xi) \dot{W}\left(
t ,\xi\right)  dt ,\\
X^{u}\left(  t ,0\right)  =X^{u}\left(  t ,1\right)  =0,\\
X^{u}\left(  t,\xi\right)  =x_{0}\left(  \xi\right)  ,
\end{array}
\right.  \label{heat equation}
\end{equation}
where $\chi_{[a,b]}$ is the indicator function of $[a,b]$ with
$0\leq a\leq b\leq 1$; $\dot{W}\left( t ,\xi\right)  $ is a
space-time white noise on $\left[ 0,T\right] \times\left[
0,1\right] $.

\noindent We introduce the cost functional
\begin{equation}
J\left(  x,u\right)  = \limsup_{T\rightarrow\infty}\dfrac{1}{T}
\mathbb{E}\int_{0}^{T}\int_{0}^{1}l\left(  \xi ,X^{u}_s\left(
\xi\right)  ,u_s(\xi)\right)  \mu\left(  d\xi\right) \, ds,
 \label{heat costo diri}
\end{equation}
where $\mu$ is a finite Borel measure on $\left[  0,1\right]  $.
An admissible control $u\left(  t  ,\xi\right)  $ is a predictable
process such that for all $t \geq 0$, and
 $\mathbb{P}$-a.s.
 $u\left(  t ,\cdot\right)
\in U:=\{v\in C\left(  \left[  0,1\right]  \right) :\left\vert
v\left( \xi\right)  \right\vert \leq\delta\}$. We denote by
$\mathcal{U}$ the set of such admissible controls. We wish to
minimize the cost  over $\mathcal{U}$, adopting the formulation of
Section \ref{optcontr}, i.e. by a change of probability in the
form of (\ref{def-ergodic-cost}). The cost introduced in
(\ref{heat costo diri}) is well defined on the space of continuous
functions on the interval $\left[ 0,1\right] $, but for an
arbitrary $\mu$\ it is not well defined on the Hilbert space of
square integrable functions.

We  suppose the following:

\begin{hypothesis}
\label{heatipotesi}
\begin{enumerate}
\item $f:\left[  0,1\right]  \times\mathbb{R} \to\mathbb{R}$ is
continuous and for every
 $\xi\in\left[0,1\right] $, $ f(\xi,\,\cdot\,)$ is decreasing.
Moreover there exist $C>0$ and $m>0$ such that for every
$\xi\in\left[0,1\right]  ,$ $x\in\mathbb{R}$,
$$ |f\left(
\xi,x\right)|\leq C(1+|x|)^m, \qquad f\left( 0,x\right)= f\left(
1,x\right)=0.
$$

\item $l:\left[  0,1\right]  \times\mathbb{R} \times
[-\delta,\delta]\rightarrow\mathbb{R}$ is continuous and bounded,
and $l(\xi,\cdot,u)$ is Lipschitz continuous uniformly with
respect to $\xi \in\left[ 0,1\right]$, $u\in [-\delta,\delta]$.

\item $x_{0}\in C\left(  \left[  0,1\right]  \right)  $,
$x_{0}(0)=x_{0}(1)=0$.
\end{enumerate}
\end{hypothesis}

\noindent To rewrite the problem in an abstract way we set
 $H=\Xi=L^{2}\left(   0,1  \right)  $
 and $E=C_0\left(\left[  0,1\right]  \right)
 =\{y\in C\left(\left[  0,1\right]  \right)\,:\, y(0)=y(1)=0\} $.
 We define an operator $A$ in $E$\ by
\[
D\left(  A\right)  =\{y\in C^{2}\left(  \left[  0,1\right]
\right)\,:\, y,y''\in C_{0}\left(  \left[  0,1\right] \right)\}
,\text{ \ \ \ \ }\left( Ay\right) \left(  \xi\right)
=\frac{\partial^{2}}{\partial\xi^{2}}y\left(  \xi\right) \text{
for  }y\in D\left(  A\right).
\]
We notice that $A$ is the generator of a $C_0$ semigroup in $E$,
admitting and extension to $H$, and $\left| e^{tA}\right|
_{L\left( E,E\right) }\leq e^{- t}$ see, for instance,  Theorem
11.3.1 in \cite{DP2}. As a consequence, $A+ F+I$ is
  dissipative in  $E$.

We set, for $x\in E$, $\xi \in [0,1]$, $z\in \Xi$, $u\in U$,
\begin{equation}
F\left(  x\right)  \left(  \xi\right)  =f\left(  \xi,x\left(
\xi\right)  \right)  ,\ \ \left(  Gz\right)  \left(  \xi\right)
 =\chi_{[a,b]}\left(  \xi\right)  z\left(
\xi\right) ,\ \ L\left(  x,u\right)
=\displaystyle\int_{0}^{1}l\left( \xi,x\left(  \xi\right) ,u\left(
\xi\right) \right)  \mu\left(  d\xi\right)  ,
\label{heatnotazioni}
\end{equation}
and let  $R$ denote the canonical imbedding of $C(  \left[
0,1\right])$ in
 $L^2(   0,1)$.

\noindent Finally $\left\{  W_{t },t \geq0\right\}  $ is a
cylindrical Wiener process in $H$ with respect to the filtration
$\left( \mathcal{F}_{t }\right) _{t \geq0}$

$ $

\noindent It is easy to verify that Hypotheses
\ref{general_hyp_forward} and \ref{hyp_W_A F(W_A)} are satisfied
(for the proof of point $4$ in Hypothesis
\ref{general_hyp_forward} and of Hypothesis \ref{hyp_W_A F(W_A)}
see again \cite{DP2} Theorem 11.3.1.).

\noindent Moreover, see for instance \cite{C}, for some $C>0$,
\[
\left| e^{tA}\right| _{L\left( H,E\right)  }\leq Ct^{-1/4}, \qquad
 t\in(0,1] ,
\]
thus Hypothesis \ref{hyp-convol-determ} holds.

\noindent Also Hypothesis \ref{Hyp-masiero} is satisfied by taking
$\Xi _{0}=\left\lbrace f\in C_0\left(  \left[  0,1\right]
\right):f(a)=f(b)=0 \right\rbrace $.

$ $ \noindent Clearly the controlled heat equation (\ref{heat
equation}) can now be written in abstract way in the Banach space
$E$ as
\begin{equation}
\left\{
\begin{array}
[c]{l}
dX_{t }^{x_0,u}=\left[  AX_{t }^{x_0,u}+F\left( X_{t
}^{x_0,u}\right) \right] dt +GRu_{t }dt +GdW_{t }\text{\ \ \ }t
\in\left[
t,T\right] \\
X^{x_0,u}_0=x_{0},
\end{array}
\right.  \label{heat eq abstract}
\end{equation}
and  the results of the previous sections can be applied to the
ergodic cost (\ref{heat costo diri}) (reformulated by a change of
probability in the form of (\ref{def-ergodic-cost})).

\noindent In particular if we define,
for all $x\in C_0([0,1])$, $z\in L^2(0,1)$, $u\in U$
(identifying $L^2(0,1)$ with its  dual)
$$\psi(x,z)=\inf _{u\in U}\left\{\int_0^1 l (\xi,x(\xi),u(\xi))
 \mu (d\xi)+ \int_a^b z(\xi) u(\xi) d\xi\right\}$$
then there exist $\overline v: E \rightarrow \mathbb{R}$
Lipschitz continuous and with $\overline v(0)=0$, $\overline \zeta : E
\rightarrow \Xi^*$ measurable and $\overline \lambda \in
\mathbb{R}$ such that if $X^{x_0}=X^{x_0,0}$ is the solution of
equation (\ref{heat eq abstract}) then $(\overline v(X^{x_0}),
\overline \zeta(X^{x_0}),\overline \lambda)$ is a solution of the
EBSDE (\ref{EBSDE}) and the characterization of the optimal
ergodic control stated in  Theorem \ref{Th-main-control} holds (and
 $\overline \lambda$ is unique in the sense of Theorem
\ref{th-uniq-lambda}).

 $ $

\noindent Moreover if $ f$ is of class $C^1(\mathbb{R})$
 (consequently $F$ will be of class ${\cal G}^1(E,E)$) and $\psi$
  is of class ${\cal G}^1(E\times \Xi^*,E)$ then by Theorem \ref{th-diff}
 $ \overline v$ is of class ${\cal G}^1(E,E)$ and, by Theorem \ref{th-EHJB},
it is a mild solution of the ergodic HJB equation (\ref{hjb}) and it holds
 $\overline \zeta=\nabla \overline v G$.

$ $

\noindent Let us then consider the particular case in which $[a,b]
=[0,1]$,  $f(x,\xi)=f(x)$ is of class $C^1$ with derivative having polynomial
growth, and satisfies $f(0)=0$,
$[f(x+h)-f(x)]h\leq - c |h|^{2+\epsilon}$
for suitable $c,\epsilon
>0$ and all $x,h\in \mathbb{R}$ (for instance, $f(x)=-x^3$).
In that case the Kolmogorov semigroup corresponding to the process
$X^{x_0}$ is strongly Feller, see
 \cite{C} and  \cite{masiero2},  and it is easy to verify that
$F$ is genuinely dissipative (see Definition \ref{gen-diss}).
Moreover we can choose $\Xi_0=C_0([0,1])$ and it turns out that
$\psi$
 is bounded
on each set $E\times B$, where $B$ is any ball of  $\Xi_0^*$. Thus
the claims of Corollaries \ref{characterization of lambda} and
\ref{boundedness of v} hold true, and in particular $\overline v$
is bounded.

$ $

\noindent Finally if we assume that $\mu$ is Lebesgue measure and
$f$ is bounded and Lipschitz we can choose
$E=\Xi=\Xi_0=H=L^2(0,1)$. Then the assumptions of Theorem
\ref{th-rec-seidler} are satisfied and we can apply Theorem
\ref{th-uniqueness} to characterize the function $\overline v$. In
particular if $f$ is of class $C^1(\mathbb{R})$ and $\psi$ is of
class ${\cal G}^1(H\times \Xi^*,H)$ then $\overline v$ is the
unique mild solution of the ergodic HJB equation (\ref{hjb}).

\end{document}